\documentclass[12pt,a4paper]{article}
\usepackage{amssymb,amsmath}
\newtheorem{theorem}{Theorem}

\usepackage[english,ngerman]{babel}

\usepackage[top=2.0cm,bottom=2.0cm,left=2.0cm,right=2.0cm]{geometry}

\begin{document}

\begin{center}
{\LARGE\bf
Series for even powers of Pi by generalization Euler's method for solving the Basel Problem
}
$$
$$
\large
{Alois Schiessl\\
Amateur Research Mathematician}
\vskip 0.5cm
{E-Mail: aloisschiessl@web.de}
\vskip 0.75cm
\selectlanguage{english}
\centerline{\today}
\end{center}
$$
$$
\selectlanguage{english}
\begin{abstract}
The purpose of this paper is to present series expansions for even powers of the number $\pi$. This is accomplished by generalizing Euler's method for solving the Basel Problem, which was published in 1735. We employ elementary symmetric polynomials, transform them into nested sums, and thereby derive nice series formulas for even powers of the number $\pi$, such as
\[
\frac{\pi^2}{3!}=\sum _{\ell_{{1}}=1}^{\infty} \frac{1}{\ell_1^2}
\;;\quad\quad
\frac{\pi^4}{5!} 
=\sum _{\ell_{{2}}=2}^{\infty} \sum _{\ell_{{1}}=1}^{\ell_{{2}}-1} 
\frac{1}{\ell_1^2 \cdot\ell_2^2}
\;;\quad\quad
\frac{\pi^6}{7!}=\sum _{\ell_{{3}}=3}^{\infty} \sum _{\ell_{{2}}=2}^{\ell_{{3}}-1}  \sum _
{\ell_{{1}}=1}^{\ell_{{2}-1}} \frac{1}{\ell_1^2 \cdot\ell_2^2 \cdot\ell_3^2} 
\;;\quad\cdots
\]
Many of these formulas do not seem to be widely known.
$$
$$
\selectlanguage{ngerman}
\qquad\qquad\qquad\qquad\qquad\qquad\qquad\qquad\textbf{Zusammenfassung}

In dieser Abhandlung stellen wir ein Verfahren vor, das die Berechnung von Reihen für geradzahlige Potenzen von $\pi$ ermöglicht. Die Grundidee ist eine Verallgemeinerung des Verfahrens von Euler, mit dem er 1735 das Basler Problem löste. Wir stellen elementar-symmetrische Polynome durch mehrfach verschachtelte Summen dar und leiten davon Reihen für geradzahlige Potenzen der Kreiszahl $\pi$ ab. Die meisten der angegebenen Reihen scheinen nicht so bekannt zu sein.
\end{abstract}

\centerline{Deutsche Version ab Seite 18}
\selectlanguage{english}
\centerline{}
\centerline{}
\centerline{***********************}
\centerline{\it In celebration of Pi Day}
\centerline{\it 3 - 14 - 2024}
\centerline{***********************}
\section{Introduction}
The Basel Problem was a highly renowned mathematical challenge in the mid-seventeenth century, initially proposed by Pietro Mengoli in 1650. This problem revolves around determining the sum of the reciprocals of the squares of positive integers. It wasn't until nearly a century later that Euler  \cite{E1}, \cite{E2} achieved success in 1735 by demonstrating that
\[
\sum_{k=1}^\infty \frac{1}{k^2}=\frac{\pi^2}{6}.
\]
Euler's original derivation of the value of $\pi^2/6$ essentially extended observations about the product representation of finite polynomials. A polynomial with degree $n$ and leading coefficient 1 can be factored as a product $\left(x-x_1\right)\ldots\left(x-x_n\right)$, where $x_1,\ldots,x_n$ are its (complex) roots. Euler made the assumption that these same properties hold true for an infinite product representation of a function. To elucidate Euler's argument, let's revisit the Taylor series expansion of the sine function.

\[ \sin(x)=x-\frac{x^3}{3!}+\frac{x^5}{5!}-+\cdots\]
Euler considered
\[
{P}\left(x\right)=\left\{
\begin{aligned}
\frac{\sin\left(x\right)}{x},\qquad &x \ne 0
\\ 
1\,\;\quad\qquad &x = 0\\ 
\end{aligned} 
\right.
\]
and interpreted $P\left(x\right)$ as a polynomial of infinite degree.
\\
Since the roots of $P\left(x\right)$ are $\{\pm\pi;\pm2\pi;\pm3\pi;\ldots\}$ he assumed that it could be factored as 
\[
\frac{\sin(x)}{x}=
\left(1-\frac{x^2}{\pi^2}\right)\left(1-\frac{x^2}{2^2\pi^2}\right)\left(1-\frac{x^2}{3^2\pi^2}\right)\ldots\]
and concluded that the following equation holds
\[
1-\frac{x^2}{3!}+\frac{x^4}{5!}-+\cdots=
\left(1-\frac{x^2}{\pi^2}\right)\left(1-\frac{x^2}{2^2\pi^2}\right)\left(1-\frac{x^2}{3^2\pi^2}\right)\ldots
\]
Euler asserted that since the left and right hand sides have the same roots and the same value at $x=0$ they are equivalent descriptions of the same function. While they are the same function, Euler’s logic does not hold. Nevertheless, Euler boldly proceeds to expand the right side product, where he considers and collects only the $x^2$ terms that emerge. This gives
\[
1-\frac{x^2}{3!}+\frac{x^4}{5!}-+...=
1-\left(\frac{1}{\pi^2}+\frac{1}{2^2\pi^2}+\frac{1}{3^2\pi^2}\right) x^2\;+-\;\ldots\]
and implies the celebrated result
\[
\frac{\pi^2}{6}=\sum_{k=1}^\infty \frac{1}{k^2}\;.
\]
A hundred years later, Karl Weierstrass finally provided a rigorous proof using techniques from complex analysis that Euler's method is indeed valid. 
$ \\ $
$ \\ $
 In this article, we generalize Euler's method to higher powers and employ nested sums to generate infinite sum representations for even powers of the number $\pi$.

\section{Nested sums}
Let $M\in\mathbb N$ and $x_{1},\ldots ,x_{M}$ be $M$ independent indeterminate. We consider the product
\begin{align*}
\left(1+x_1 t\right)\cdots\left(1+x_M t\right)=\prod_{k=1}^{M}(1+x_{k}t)
\end{align*}
By expanding it, collecting and rearranging the terms, we have the identity
\begin{align*}
\prod_{k=1}^{M}(1+x_{k}t)
=1+\sigma_{M,1}t+\sigma_{M,2}t^2+\sigma_{M,3}t^3+\ldots+\sigma_{M,M}t^M=1+\sum_{k=1}^{M}\sigma_{M,k}\cdot t^{k}
\end{align*}
with the elementary symmetric polynomials
\begin{align*}
\sigma_{M,k}&=\begin{cases}
{\sum\limits_{1 \le {\ell_1} <\cdots < {\ell_k} \le M} {x_{{\ell_1}}\cdots x_{{\ell_k}}} } & (1\leq k\leq M) \\
\\
\qquad 0 & (M<k)
\end{cases}
\end{align*}
respectively.
Of course, we can write them down explicitly. If $k\in\mathbb \{1,2,3,...,M\}$ then the elementary symmetric polynomials are
\[
\sigma_{M,1}=\sum\limits_{1 \le {\ell_1} \le M} {x{}_{{\ell_1}}}
\]
\[
\sigma_{M,2}=\sum\limits_{1 \le {\ell_1} < {\ell_2} \le M} 
x_{\ell_1}\cdot x_{\ell_2}
\]
\[
\sigma_{M,3}=\sum\limits_{1 \le {\ell_1} < {\ell_2} < {\ell_3} \le M} x_{\ell_1}\cdot x_{\ell_2}\cdot x_{\ell_3}
\]
\[
\centerline{\vdots}
\]
\[
\sigma_{M,k}={\sum\limits_{1 \le {\ell_1} <\cdots < {\ell_k} \le N} {x_{{\ell_1}}\cdots x_{{\ell_k}}} }
\]

\[
\centerline{\vdots}
\]
\[
\sigma_{M,M}={\sum\limits_{1 \le {\ell_1} <\cdots < {\ell_M} \le M} {x_{{\ell_1}}\cdots x_{{\ell_M}}} }=\prod_{\ell_k=1}^{M}x_{\ell_k}
\]
\\
The following lemma \cite{EF} regarding the properties of the elementary symmetric polynomials turns out to be highly beneficial:
$ \\ $
$ \\ $
\textbf{Lemma 1}
\begin{align*}
\sigma_{M,k}&=0, \;M < k
\\
\sigma_{M,0}&=1, \;\; 0 \;\leq M
\\
\sigma_{M+1,k}&=\sigma_{M,k}+x_{M+1}\; \sigma_{M,k-1},\;\;0 < k \le M+1
\end{align*}
\\
The elementary symmetric polynomial in its general form
\[
{\sum\limits_{1 \le {\ell_1} <\cdots < {\ell_k} \le M} {x_{{\ell_1}}\cdots x_{{\ell_k}}} }
\]
cannot be easily utilized for computational calculation. Therefore it is desirable to transform them into an appropriate form. Although this may seem like a challenging task, it can be accomplished step by step using nested sums. We only need to break down the multi-inequality 
\[
1 \le {\ell_1} < {\ell_2} <\cdots < {\ell_{k-1}} < {\ell_k} \le M
\]
into $k$ single inequalities and looking for the starting and ending index.

$ \\ $
Starting with $k=1$, we get from $1 \le {\ell_1} \le M$ immediately $1 \le\ell_1$ and $\ell_1 \le M$, that we obtain the simple sum
\[
\sigma_{M,1}
=\sum_{\ell_1=1}^{M}x_{\ell_1}
=x_1+\ldots+x_M
\]
\\
Having $k=2$, we split the chained notation $1 \le\ell_1 <\ell_2 \le M$ into two inequalities
\begin{align*}
1 \le\ell_1 &\;\;\text{and}\;\;\ell_1 \le\ell_2-1\\
2 \le\ell_2 &\;\;\text{and}\;\;\ell_2 \le M
\end{align*}
So, we obtain the double sum
\[
\sigma_{M,2}
=\underbrace {\sum _{\ell_{{2}}=2}^{M} \left( \sum _{\ell_{{1}}=1}^{\ell_{{2}}-1}x_{{\ell_{{1}}}
}\cdot x_{{\ell_{{2}}}} \right)}_{2\;nested\;sums}\\
=
x_1x_2+\ldots+x_{M-1}x_M
\]
\\
For $k=3$, we generate from the chained notation $1 \le\ell_1 <\ell_2 <\ell_3 \le M$ three inequalities
\begin{align*}
1 \le\ell_1 &\;\;\text{and}\;\;\ell_1 \le\ell_2-1\\
2 \le\ell_2 &\;\;\text{and}\;\;\ell_2 \le\ell_3-1\\
3 \le\ell_3 &\;\;\text{and}\;\;\ell_3 \le M
\end{align*}
that yields the triple sum
\[
\sigma_{M,3}
=\underbrace {\sum _{\ell_{{3}}=3}^{M} \left( \sum _{\ell_{{2}}=2}^{\ell_{{3}}-1} \left( 
\sum _{\ell_{{1}}=1}^{\ell_{{2}}-1} x_{{\ell_{{1}}}} \cdot x_{{\ell_{{2}}}}\cdot x_{{\ell_{{3}}}} \right)  \right)}_{3\;nested\;sums}
\]
\[
=x_1\cdot x_2\cdot x_3+\cdots+x_{M-2}\cdot x_{M-1}\cdot x_M
\]
\\
In the same way for $k=4$, we generate from the chained notation $1 \le\ell_1 <\ell_2 <\ell_3 <\ell_4 \le M$ four inequalities
\begin{align*}
1 \le\ell_1 &\;\;\text{and}\;\;\ell_1 \le\ell_2-1\\
2 \le\ell_2 &\;\;\text{and}\;\;\ell_2 \le\ell_3-1\\
3 \le\ell_3 &\;\;\text{and}\;\;\ell_3 \le\ell_4-1\\
4 \le\ell_4 &\;\;\text{and}\;\;\ell_4 \le M
\end{align*}
and obtain the quadruple sum
\[\sigma_{M,4}
=\underbrace{\sum _{\ell_{{4}}=4}^{M} \left(\sum _{\ell_{{3}}=3}^{\ell_4-1} \left( \sum _{\ell_{{2}}=2}^{\ell_{{3}}-1} \left( 
\sum _{\ell_{{1}}=1}^{\ell_{{2}}-1} x_{{\ell_{{1}}}} \cdot x_{{\ell_{{2}}}}\cdot x_{{\ell_3}}\cdot x_{{\ell_4}} \right)  \right)\right)}_{4\;nested\;sums}
\]
\[
=x_1\cdot x_2\cdot x_3\cdot x_4+\cdots+x_{M-3}\cdot x_{M-2}\cdot x_{M-1}\cdot x_M
\]
More generally, for $k=M$ the chained notation 
\[1 \le {\ell_1} <\ell_2 \cdots <l_{k-1}  < {\ell_k} \le M\]
leads to $k$ inequalities 
\begin{align*}
1 \le\ell_1 &\;\;\text{and}\;\;\ell_1 \le\ell_2-1\\
2 \le\ell_2 &\;\;\text{and}\;\;\ell_2 \le\ell_3-1\\
&\quad\,\vdots\\
k-1 \le\ell_{k-1} &\;\;\text{and}\;\;\ell_{k-1} \le\ell_k-1\\
k \le\ell_k &\;\;\text{and}\;\;\ell_k \le M
\end{align*}
From these inequalities we can express the following $k$ nested sum:
\[
\sigma_{M,k}
=\underbrace {\sum\limits_{{\ell_k} = k}^M {\left( {\sum\limits_{{\ell_{k - 1}} = k - 1}^{{\ell_k} - 1}  \cdots  \left( {\sum\limits_{{\ell_2} = 2}^{{\ell_3} - 1} {\left( {\sum\limits_{{\ell_1} = 1}^{{\ell_{2} - 1}} x_{\ell_1}\cdot x_{\ell_2} \cdots x_{\ell_{k-1}} \cdot x_{\ell_k} } \right)} } \right) \cdots } \right)} }_{k\;nested\;sums}
\]
\[={x_1\cdots x_k}+\cdots+x_{M-k+1}\cdots x_M 
\]\\
Note: the number of terms
of $\sigma_{M,k}$ is given by $\#\sigma_{M,k}=\binom{M}{k}.$
$ \\ $
$ \\ $
We now present the following theorem:
\begin{theorem}
Let $M\in\mathbb N$ and $x_{1},\ldots ,x_{M}$ be $M$ independent indeterminate. The product
\begin{align*}
\prod_{k=1}^{M}(1+x_{k}t)=\left(1+x_1 t\right)\cdots\left(1+x_M t\right)
\end{align*}
can be expressed by a polynomial with coefficients consisting of multi sums
\begin{align*}
\prod_{k=1}^{M}(1+x_{k}t)=1
+&\left[\sum_{\ell_1=1}^{M}x{_{\ell_1}}\right] t
\\
+&\left[\sum _{\ell_{{2}}=2}^{M} \left( \sum _{\ell_{{1}}=1}^{\ell_{{2}}-1}x_{{\ell_{{1}}}
} x_{{\ell_{{2}}}} \right)\right] t^2
\\
+&\left[\sum _{\ell_{{3}}=3}^{M} \left( \sum _{\ell_{{2}}=2}^{\ell_{{3}}-1} \left( 
\sum _{\ell_{{1}}=1}^{\ell_{{2}}-1}x_{{\ell_{{1}}}} x_{{\ell_{{2}}}} x_{{\ell_{{3}}}} \right)  \right)\right] t^3
\\
+&\left[{\sum _{\ell_{{4}}=4}^{M} \left(\sum _{\ell_{{3}}=3}^{\ell_4-1} \left( \sum _{\ell_{{2}}=2}^{\ell_{{3}}-1} \left( 
\sum _{\ell_{{1}}=1}^{\ell_{{2}}-1} x_{{\ell_{{1}}}}  x_{{\ell_{{2}}}} x_{{\ell_3}} x_{{\ell_4}} \right)  \right)\right)}\right] t^4
\\
&\qquad\qquad\qquad\qquad{\vdots}
\\
+&\left[\sum\limits_{{\ell_M} = M}^M {\left( {\sum\limits_{{\ell_{M - 1}} = M - 1}^{{\ell_M} - 1}  \cdots \left( {\sum\limits_{{\ell_2} = 2}^{{\ell_3} - 1} {\left( {\sum\limits_{{\ell_1} = 1}^{{\ell_{2} - 1}} x_{\ell_1} x_{\ell_2} \cdots x_{\ell_{M-1}} x_{\ell_M} } \right)} } \right) \cdots } \right)}\right]t^M
\end{align*}
\end{theorem}
The proof of the theorem is accomplished by an induction argument on $M$. While not particularly challenging, some effort is required. We defer the proof to the Appendix.

\section{From polynomial to sinc function}
We consider the Polynomial
\begin{align*}
\prod_{k=1}^{M}(1+x_{k}t)=1
+&\left[\sum_{\ell_1=1}^{M}x{_{\ell_1}}\right] t
\\
+&\left[\sum _{\ell_{{2}}=2}^{M} \left( \sum _{\ell_{{1}}=1}^{\ell_{{2}}-1}x_{{\ell_{{1}}}
} x_{{\ell_{{2}}}} \right)\right] t^2
\\
+&\left[\sum _{\ell_{{3}}=3}^{M} \left( \sum _{\ell_{{2}}=2}^{\ell_{{3}}-1} \left( 
\sum _{\ell_{{1}}=1}^{\ell_{{2}}-1}x_{{\ell_{{1}}}} x_{{\ell_{{2}}}} x_{{\ell_{{3}}}} \right)  \right)\right] t^3
\\
+&\left[{\sum _{\ell_{{4}}=4}^{M} \left(\sum _{\ell_{{3}}=3}^{\ell_4-1} \left( \sum _{\ell_{{2}}=2}^{\ell_{{3}}-1} \left( 
\sum _{\ell_{{1}}=1}^{\ell_{{2}}-1} x_{{\ell_{{1}}}}  x_{{\ell_{{2}}}} x_{{\ell_3}} x_{{\ell_4}} \right)  \right)\right)}\right] t^4
\\
&\qquad\qquad\qquad\qquad{\vdots}
\\
+&\left[\sum\limits_{{\ell_M} = M}^M {\left( {\sum\limits_{{\ell_{M - 1}} = M - 1}^{{\ell_M} - 1}  \cdots \left( {\sum\limits_{{\ell_2} = 2}^{{\ell_3} - 1} {\left( {\sum\limits_{{\ell_1} = 1}^{{\ell_{2} - 1}} x_{\ell_1} x_{\ell_2} \cdots x_{\ell_{M-1}} x_{\ell_M} } \right)} } \right) \cdots } \right)}\right]t^M
\end{align*}
On the left hand side we substitute $t=-x^2$ and $x_k=\frac{1}{k^2}$. That means on the right hand we have to substitute also $t=-x^2$ and $x_{\ell_k}=\frac{1}{\ell_k^2}$. That results in
\begin{align*}
\prod_{k=1}^{M}(1-\frac{x^2}{k^2})=1-&\left[\sum_{\ell_1=1}^{N}\frac{1}{\ell_1^2} \right] x^2\\
+&\left[\sum _{\ell_{{2}}=2}^{N} \left( \sum _{\ell_{{1}}=1}^{\ell_{{2}}-1}\frac{1}{\ell_1^2\cdot \ell_2^2} \right)\right] x^4\\
-&\left[\sum _{\ell_{{3}}=3}^{N} \left( \sum _{\ell_{{2}}=2}^{\ell_{{3}}-1} \left( 
\sum _{\ell_{{1}}=1}^{\ell_{{2}}-1}
\frac{1}{\ell_1^2\cdot \ell_2^2\cdot \ell_3^2} \right)  \right)\right] x^6\\
+&\left[\sum _{\ell_{{4}}=4}^{N} \left( \sum _{\ell_{{3}}=3}^{\ell_{{4}}-1}
 \left( \sum _{\ell_{{2}}=2}^{\ell_{{3}}-1} \left( \sum _{\ell_{{1}}=1}^{\ell_{{2}
}-1}
\frac{1}{\ell_1^2\cdot \ell_2^2\cdot \ell_3^2\cdot \ell_4^2} \right) 
 \right)  \right)\right] x^8
\\
&\qquad\qquad\qquad\qquad{\vdots}
\\
+\left(-1\right)^M
&\left[\sum\limits_{{\ell_M} = M}^M {\left( {\sum\limits_{{\ell_{M - 1}} = M - 1}^{{\ell_M} - 1}  \cdots  \left( {\sum\limits_{{\ell_2} = 2}^{{\ell_3} - 1} {\left( {\sum\limits_{{\ell_1} = 1}^{{\ell_{2 - 1}}} \frac{1}{\ell_1^2\cdot \ell_2^2 \cdots\ell_{M-1}^2\cdot \ell_M^2} } \right)} } \right) \cdots } \right)}\right]x^{2M}
\end{align*}
\\
Now, let us investigate what happens when taking the limit $M \to \infty $ on both sides.
$$
$$
The left hand side
\[
\underset{M\to \infty }{\mathop{\lim }}\,\prod_{k=1}^{M}(1-\frac{x^2}{k^2})
=\prod_{k=1}^{\infty}(1-\frac{x^2}{k^2})
\]
results in the infinite product formula for the sinc function.
The $sinc$ function is defined as 
\[
\text{sinc}\left( \pi x\right)
=\left\{
\begin{aligned}
\frac{\sin\left( \pi x \right)}{\pi x}\,,\qquad &x \ne 0
\\ 
1\;,\;\quad\qquad &x = 0 
\end{aligned} 
\right.
\]
with power series expansion
\[
\frac{\sin\left(\pi x\right)}{\pi x}=1-\frac{\pi^2}{3!}x^2+\frac{\pi^4}{5!}x^4-\frac{\pi^6}{7!}x^6+-\cdots
\]
This can be easily derived from the Weierstrass Factorization Theorem for the $\sin$ function.
$$
$$
On the right hand side when taking the limit $M \to \infty $, we obtain a power series with coefficients consisting of infinite nested sums
\begin{align*}
1-&\left[\sum_{\ell_1=1}^{\infty}\frac{1}{\ell_1^2} \right] x^2\\
+&\left[\sum _{\ell_{{2}}=2}^{\infty} \left( \sum _{\ell_{{1}}=1}^{\ell_{{2}}-1}\frac{1}{\ell_1^2\cdot \ell_2^2} \right)\right] x^4\\
-&\left[\sum _{\ell_{{3}}=3}^{\infty} \left( \sum _{\ell_{{2}}=2}^{\ell_{{3}}-1} \left( 
\sum _{\ell_{{1}}=1}^{\ell_{{2}}-1}
\frac{1}{\ell_1^2\cdot \ell_2^2\cdot \ell_3^2} \right)  \right)\right] x^6\\
+&\left[\sum _{\ell_{{4}}=4}^{\infty} \left( \sum _{\ell_{{3}}=3}^{\ell_{{4}}-1}
 \left( \sum _{\ell_{{2}}=2}^{\ell_{{3}}-1} \left( \sum _{\ell_{{1}}=1}^{\ell_{{2}
}-1}
\frac{1}{\ell_1^2\cdot \ell_2^2\cdot \ell_3^2\cdot \ell_4^2} \right) 
 \right)  \right)\right] x^8-+\ldots
\end{align*}
\\
Combining both sides yields a new series expansion of the sinc function 
\begin{align*}
\frac{\sin\left(\pi x\right)}{\pi x}=1-&\left[\sum_{\ell_1=1}^{\infty}\frac{1}{\ell_1^2} \right] x^2\\
+&\left[\sum _{\ell_{{2}}=2}^{\infty} \left( \sum _{\ell_{{1}}=1}^{\ell_{{2}}-1}\frac{1}{\ell_1^2\cdot \ell_2^2} \right)\right] x^4\\
-&\left[\sum _{\ell_{{3}}=3}^{\infty} \left( \sum _{\ell_{{2}}=2}^{\ell_{{3}}-1} \left( 
\sum _{\ell_{{1}}=1}^{\ell_{{2}}-1}
\frac{1}{\ell_1^2\cdot \ell_2^2\cdot \ell_3^2} \right)  \right)\right] x^6\\
+&\left[\sum _{\ell_{{4}}=4}^{\infty} \left( \sum _{\ell_{{3}}=3}^{\ell_{{4}}-1}
 \left( \sum _{\ell_{{2}}=2}^{\ell_{{3}}-1} \left( \sum _{\ell_{{1}}=1}^{\ell_{{2}
}-1}
\frac{1}{\ell_1^2\cdot \ell_2^2\cdot \ell_3^2\cdot \ell_4^2} \right) 
 \right)  \right)\right] x^8-+\ldots
\end{align*}
\section{Powers of number Pi}
Plugging in the power series expansion of the sinc function into the left-hand side of the above equation yields
\begin{align*}
1-&\frac{\pi^2}{3!}x^2+\frac{\pi^4}{5!}x^4-\frac{\pi^6}{7!}x^6+-\cdots
\\
\\
=1-&\left[\sum_{\ell_1=1}^{\infty}\frac{1}{\ell_1^2} \right] x^2
+\left[\sum _{\ell_{{2}}=2}^{\infty} \left( \sum _{\ell_{{1}}=1}^{\ell_{{2}}-1}\frac{1}{\ell_1^2\cdot \ell_2^2} \right)\right] x^4
-\left[\sum _{\ell_{{3}}=3}^{\infty} \left( \sum _{\ell_{{2}}=2}^{\ell_{{3}}-1} \left( 
\sum _{\ell_{{1}}=1}^{\ell_{{2}}-1}
\frac{1}{\ell_1^2\cdot \ell_2^2\cdot \ell_3^2} \right)  \right)\right] x^6\\
+&\left[\sum _{\ell_{{4}}=4}^{\infty} \left( \sum _{\ell_{{3}}=3}^{\ell_{{4}}-1}
 \left( \sum _{\ell_{{2}}=2}^{\ell_{{3}}-1} \left( \sum _{\ell_{{1}}=1}^{\ell_{{2}
}-1}
\frac{1}{\ell_1^2\cdot \ell_2^2\cdot \ell_3^2\cdot \ell_4^2} \right) 
 \right)  \right)\right] x^8-+\ldots
\end{align*}
By the uniqueness theorem for power series, it follows that the two series coincide. If two power series are equal, then the corresponding coefficients must be equal.
$ \\ $
$ \\ $
On equating the coefficients of power $x^2$ we obtain
\[
\frac{\pi^2}{3!}=\sum _{\ell_{{1}}=1}^{\infty} \frac{1}{\ell_1^2} 
\]
This is the celebrated result due to Euler, published in 1734. It`s intriguing why Euler did not investigating higher powers. Probably because he was solely interested in series of the form
\[
\sum_{k=1}^\infty \frac{1}{k^n}
\quad\text{with}\; n=2,3,\dots
\]
We will now proceed and explore higher powers.
\\
On equating the coefficients of power $x^4$, we find on the left side a rational multiple of $\pi^4$ and on the right hand side a double sum
\[
\frac{\pi^4}{5!} 
=\sum _{\ell_{{2}}=2}^{\infty} \left( \sum _{\ell_{{1}}=1}^{\ell_{{2}}-1} 
\frac{1}{\ell_1^2 \cdot\ell_2^2}\right)
\]
We define the right hand side limit 
\[
\sum _{\ell_{{2}}=2}^{\infty} \left( \sum _{\ell_{{1}}=1}^{\ell_{{2}}-1} 
\frac{1}{\ell_1^2 \cdot\ell_2^2}\right)
\]
as the limit of the partial sums
\[
\sum _{\ell_{{2}}=2}^{2} \left( \sum _{\ell_{{1}}=1}^{\ell_{{2}}-1}{\frac {1}{{\ell_{{1}}}^{2}{\ell_{{2}}}^{2}}} \right) 
=
\frac{1}{1^2\cdot 2^2}
\]
\[
\sum _{\ell_{{2}}=2}^{3} \left( \sum _{\ell_{{1}}=1}^{\ell_{{2}}-1}{\frac {1}{{\ell_{{1}}}^{2}{\ell_{{2}}}^{2}}} \right) 
=\frac{1}{1^2\cdot 2^2}+\frac{1}{1^2\cdot 3^2}+\frac{1}{2^2\cdot 3^2}
\]
\[
\sum _{\ell_{{2}}=2}^{4} \left( \sum _{\ell_{{1}}=1}^{\ell_{{2}}-1}{\frac {1}{{\ell_{{1}}}^{2}{\ell_{{2}}}^{2}}} \right) 
=
\frac{1}{1^2\cdot 2^2}+\frac{1}{1^2\cdot 3^2}+\frac{1}{1^2\cdot 4^2}+
\frac{1}{2^2\cdot 3^2}+\frac{1}{2^2\cdot 4^2}+
\frac{1}{3^2\cdot 4^2} 
\]
\[
\centerline{\vdots}
\]
Next, equating the coefficients of power $x^6$ gives a rational multiple of $\pi^6$ and a triple sum
\[
\frac{\pi^6}{7!}=\sum _{\ell_{{3}}=3}^{\infty} \left( \sum _{\ell_{{2}}=2}^{\ell_{{3}}-1} \left( \sum _
{\ell_{{1}}=1}^{\ell_{{2}-1}} \frac{1}{\ell_1^2 \cdot\ell_2^2 \cdot\ell_3^2}  \right)  \right) 
\]
We also define the limit
\[
\sum _{\ell_{{3}}=3}^{\infty} \left( \sum _{\ell_{{2}}=2}^{\ell_{{3}}-1} \left( \sum _
{\ell_{{1}}=1}^{\ell_{{2}-1}} \frac{1}{\ell_1^2 \cdot\ell_2^2 \cdot\ell_3^2}  \right)  \right)
\]
as the limit of the partial sums
\[
\sum _{\ell_{{3}}=3}^{3}
\left( \sum _{\ell_{{2}}=2}^{\ell_{{3}}-1}
\left( \sum _{\ell_{{1}}=1}^{\ell_{{2}}-1}{\frac {1}{{\ell_{{1}}}^{2}{\ell_{{2}}}^{2}{\ell_{{3}}}^{2}}} \right)  \right)
=\frac{1}{1^2 \cdot 2^2 \cdot 3^2}
\]
\begin{align*}
\sum _{\ell_{{3}}=3}^{4}
\left( \sum _{\ell_{{2}}=2}^{\ell_{{3}}-1}
\left( \sum _{\ell_{{1}}=1}^{\ell_{{2}}-1}{\frac {1}{{\ell_{{1}}}^{2}{\ell_{{2}}}^{2}{\ell_{{3}}}^{2}}} \right)  \right)
&=
 \frac{1}{1^2 \cdot 2^2 \cdot 3^2}
+\frac{1}{1^2 \cdot 2^2 \cdot 4^2}
+\frac{1}{1^2 \cdot 3^2 \cdot 4^2}
+\frac{1}{2^2 \cdot 3^2 \cdot 4^2}
\end{align*}
\begin{align*}
\sum _{\ell_{{3}}=3}^{5}
\left( \sum _{\ell_{{2}}=2}^{\ell_{{3}}-1}
\left( \sum _{\ell_{{1}}=1}^{\ell_{{2}}-1}{\frac {1}{{\ell_{{1}}}^{2}{\ell_{{2}}}^{2}{\ell_{{3}}}^{2}}} \right)\right)
&=
 \frac{1}{1^2 \cdot 2^2 \cdot 3^2}
+\frac{1}{1^2 \cdot 2^2 \cdot 4^2}
+\frac{1}{1^2 \cdot 2^2 \cdot 5^2}
+\frac{1}{1^2 \cdot 3^2 \cdot 4^2}
+\frac{1}{1^2 \cdot 3^2 \cdot 5^2}\\
&+\frac{1}{1^2 \cdot 4^2 \cdot 5^2}
+\frac{1}{2^2 \cdot 3^2 \cdot 4^2}
+\frac{1}{2^2 \cdot 3^2 \cdot 5^2}
+\frac{1}{2^2 \cdot 4^2 \cdot 5^2}
+\frac{1}{3^2 \cdot 4^2 \cdot 5^2}
\end{align*}
\[\centerline{\vdots}\]

Continuing with equating coefficients of powers yields
\[
\frac{{{\pi ^8}}}{{9!}} = \;\underbrace {\sum\limits_{{\ell _4} = 4}^{\infty} {\left( {\sum\limits_{{\ell _3} = 3}^{{\ell _4} - 1} {\left( {\sum\limits_{{\ell _2} = 2}^{{\ell _3} - 1} {\left( {\sum\limits_{{\ell _1} = 1}^{{\ell _{2 - 1}}} {\frac{1}{{\ell _1^2 \cdot \ell _2^2 \cdot \ell _3^2 \cdot \ell _4^2}}} } \right)} } \right)} } \right)}}_{4\;nested\;sums}
\]
\[
\frac{\pi^{10}}{11!}=\;\underbrace {\sum _{\ell_{{5}}=5}^{\infty}\left(\sum _{\ell_{{4}}=4}^{\ell_{{5}}-1}\left(\sum _{\ell_{{3}}=3}^{\ell_{{4}}-1} \left( \sum _{\ell_{{2}}=2}^{\ell_{{3}}-1} \left( \sum_
{\ell_{{1}}=1}^{\ell_{{2}-1}} \frac{1}{\ell_1^2 \cdot\ell_2^2 \cdot\ell_3^2 \cdot\ell_4^2 \cdot\ell_5^2}  \right)  \right) \right) \right) }_{5\;nested\;sums}
\]
\[
\frac{\pi^{12}}{13!}=\;\underbrace {\sum _{\ell_{{6}}=6}^{\infty}\left(\sum _{\ell_{{5}}=5}^{\ell_{{6}}-1}
\left(\sum _{\ell_{{4}}=4}^{\ell_{{5}}-1}
\left(\sum _{\ell_{{3}}=3}^{\ell_{{4}}-1}
\left(\sum _{\ell_{{2}}=2}^{\ell_{{3}}-1}
\left( \sum_{\ell_{{1}}=1}^{\ell_{{2}-1}}
\frac{1}{\ell_1^2 \cdot\ell_2^2 \cdot\ell_3^2 \cdot\ell_4^2 \cdot\ell_5^2\cdot\ell_6^2}
\right)
\right)
\right)
\right)
\right)}_{6\;nested\;sums}
\]
\[
\frac{\pi^{14}}{15!}=\;\underbrace {\sum _{\ell_{{7}}=7}^{\infty}\left(\sum _{\ell_{{6}}=6}^{\ell_{{7}}-1}
\left(\sum _{\ell_{{5}}=5}^{\ell_{{6}}-1}
\left(\sum _{\ell_{{4}}=4}^{\ell_{{5}}-1}
\left(\sum _{\ell_{{3}}=3}^{\ell_{{4}}-1}
\left(\sum _{\ell_{{2}}=2}^{\ell_{{3}}-1}
\left( \sum_{\ell_{{1}}=1}^{\ell_{{2}-1}}
\frac{1}{\ell_1^2 \cdot\ell_2^2 \cdot\ell_3^2 \cdot\ell_4^2 \cdot\ell_5^2\cdot\ell_6^2}
\right)
\right)
\right)
\right)
\right)
\right)}_{7\;nested\;sums}
\]
\[
\vdots
\]
In the general case $n\in\mathbb N$, we obtain
\[
\frac{{{\pi ^{2n}}}}{{\left( {2n + 1} \right)!}} =  \;\underbrace {\sum\limits_{{\ell_n} = n}^{\infty} {\left( {\sum\limits_{{\ell_{n - 1}} = n - 1}^{{\ell_n} - 1}  \cdots  \left( {\sum\limits_{{\ell_2} = 2}^{{\ell_3} - 1} {\left( {\sum\limits_{{\ell_1} = 1}^{{\ell_{2 - 1}}} {\frac{1}{{\ell_1^2 \cdots\ell_n^2}}} } \right)} } \right) \cdots } \right)} }_{n\;nested\;sums}
\]
The limit of the sum
\[
\sum\limits_{{\ell_n} = n}^{\infty} {\left( {\sum\limits_{{\ell_{n - 1}} = n - 1}^{{\ell_n} - 1}  \cdots  \left( {\sum\limits_{{\ell_2} = 2}^{{\ell_3} - 1} {\left( {\sum\limits_{{\ell_1} = 1}^{{\ell_{2 - 1}}} {\frac{1}{{\ell_1^2 \cdots\ell_n^2}}} } \right)} } \right) \cdots } \right)} 
\]
\\
is defined as the limit of the partial sums
\[
\sum\limits_{{\ell_n} = n}^{n} {\left( {\sum\limits_{{\ell_{n - 1}} = n - 1}^{{\ell_n} - 1}  \cdots  \left( {\sum\limits_{{\ell_2} = 2}^{{\ell_3} - 1} {\left( {\sum\limits_{{\ell_1} = 1}^{{\ell_{2 - 1}}} {\frac{1}{{\ell_1^2 \cdots\ell_n^2}}} } \right)} } \right) \cdots } \right)} 
\]
$$
$$
\[
\sum\limits_{{\ell_n} = n}^{n+1} {\left( {\sum\limits_{{\ell_{n - 1}} = n - 1}^{{\ell_n} - 1}  \cdots  \left( {\sum\limits_{{\ell_2} = 2}^{{\ell_3} - 1} {\left( {\sum\limits_{{\ell_1} = 1}^{{\ell_{2 - 1}}} {\frac{1}{{\ell_1^2 \cdots\ell_n^2}}} } \right)} } \right) \cdots } \right)} 
\]
$$
$$
\[
\sum\limits_{{\ell_n} = n}^{n+2} {\left( {\sum\limits_{{\ell_{n - 1}} = n - 1}^{{\ell_n} - 1}  \cdots  \left( {\sum\limits_{{\ell_2} = 2}^{{\ell_3} - 1} {\left( {\sum\limits_{{\ell_1} = 1}^{{\ell_{2 - 1}}} {\frac{1}{{\ell_1^2 \cdots\ell_n^2}}} } \right)} } \right) \cdots } \right)} 
\]
\[
\vdots
\]
Now let's focus on the convergence of the sums. Since the sums only consists of positive terms, it is easily verified by comparison test.
$ \\ $
$ \\ $
The convergence of the sum
\[
\sum _{\ell_{{1}}=1}^{\infty} \frac{1}{\ell_1^2} 
\]
is well known and requires no further action.
$ \\ $
$ \\ $
For higher sums we use the well known sequence of series \cite{KK}
\[
\sum _{\ell_{{2}}=1}^{\infty} \left( \sum _{\ell_{{1}}=1}^{\infty} 
\frac{1}{\ell_1^2 \cdot\ell_2^2}\right)=\left(\frac{\pi^2}{6}\right)^2;
\]
\[
\sum _{\ell_{{3}}=1}^{\infty} \left( \sum _{\ell_{{2}}=1}^{\infty} \left( \sum _
{\ell_{{1}}=1}^{\infty} \frac{1}{\ell_1^2 \cdot\ell_2^2 \cdot\ell_3^2}  \right)  \right)=\left(\frac{\pi^2}{6}\right)^3;
\]

\[
\vdots
\]
\[
\sum\limits_{{\ell_n} = 1}^{\infty} {\left( {\sum\limits_{{\ell_{n - 1}} =1}^{\infty}  \cdots  \left( {\sum\limits_{{\ell_2} = 1}^{\infty} {\left( {\sum\limits_{{\ell_1} = 1}^{\infty} {\frac{1}{{\ell_1^2 \cdots\ell_n^2}}} } \right)} } \right) \cdots } \right)}
=\left(\frac{\pi^2}{6}\right)^n\]
\\
Comparing it with the derived series, we get
\[
\sum _{\ell_{{2}}=2}^{\infty} \left( \sum _{\ell_{{1}}=1}^{\ell_{{2}}-1} 
\frac{1}{\ell_1^2 \cdot\ell_2^2}\right)
<\sum _{\ell_{{2}}=1}^{\infty} \left( \sum _{\ell_{{1}}=1}^{\infty} 
\frac{1}{\ell_1^2 \cdot\ell_2^2}\right)=\left(\frac{\pi^2}{6}\right)^2;
\]
\[
\sum _{\ell_{{3}}=3}^{\infty} \left( \sum _{\ell_{{2}}=2}^{\ell_{{3}}-1} \left( \sum _
{\ell_{{1}}=1}^{\ell_{{2}-1}} \frac{1}{\ell_1^2 \cdot\ell_2^2 \cdot\ell_3^2}  \right)  \right) 
<
\sum _{\ell_{{3}}=1}^{\infty} \left( \sum _{\ell_{{2}}=1}^{\infty} \left( \sum _
{\ell_{{1}}=1}^{\infty} \frac{1}{\ell_1^2 \cdot\ell_2^2 \cdot\ell_3^2}  \right)  \right)=\left(\frac{\pi^2}{6}\right)^3;
\]
\[
\vdots
\]
More generally, for $n\in\mathbb N$, we obtain
\begin{align*}
&\sum\limits_{{\ell_n} = n}^{\infty} {\left( {\sum\limits_{{\ell_{n - 1}} = n - 1}^{{\ell_n} - 1}  \cdots  \left( {\sum\limits_{{\ell_2} = 2}^{{\ell_3} - 1} {\left( {\sum\limits_{{\ell_1} = 1}^{{\ell_{2 - 1}}} {\frac{1}{{\ell_1^2 \cdots\ell_n^2}}} } \right)} } \right) \cdots } \right)}\\
<
&\sum\limits_{{\ell_n} = 1}^{\infty} {\left( {\sum\limits_{{\ell_{n - 1}} =1}^{\infty}  \cdots  \left( {\sum\limits_{{\ell_2} = 1}^{\infty} {\left( {\sum\limits_{{\ell_1} = 1}^{\infty} {\frac{1}{{\ell_1^2 \cdots\ell_n^2}}} } \right)} } \right) \cdots } \right)}
=\left(\frac{\pi^2}{6}\right)^n
\end{align*}
We observe, that every $n$-nested sum is bounded by the value $\left(\frac{\pi^2}{6}\right)^n $, which implies convergence.

\section{Computational results with MAPLE}
The calculations are performed using MAPLE Release 2024.0. MAPLE is a computational environment that supports both symbolic and numeric computation, as well as serving as a multi-paradigm programming language. \\
If you have access to a MAPLE environment, you can directly copy and paste the provided code into a Maple session. 
$ \\ $
$ \\ $
We start with the sum
\[
\sum _{\ell_{{1}}=1}^{\infty} \frac{1}{\ell_1^2} 
\]
and use the MAPLE Code

\begin{verbatim}
sum(1/l[1]^2,l[1]=1..infinity);
\end{verbatim}
Copying and pasting it in a MAPLE session yields the desired result
\[
\frac{1}{6}\,{\pi }^{2}
\]
Similarly, we perform the same process for the double sum
\[
\sum _{\ell_{{2}}=2}^{\infty} \left( \sum _{\ell_{{1}}=1}^{\ell_{{2}}-1} 
\frac{1}{\ell_1^2 \cdot\ell_2^2}\right)
\]
The MAPLE code for the double sum is as follows: 

\begin{verbatim}
sum(sum(1/(l[1]^2*l[2]^2),
l[1]=1..l[2]-1),l[2]=2..infinity);
\end{verbatim}
MAPLE fails to provide the expected results $\frac{\pi^4}{5!}.$ Instead, it gives an expression using the $\Psi$-function:
\[
\sum _{l_{{2}}=2}^{\infty }-{\frac {\Psi \left( 1,l_{{2}} \right) }{{l
_{{2}}}^{2}}}+\frac{1}{6}\,{\frac {{\pi }^{2}}{{l_{{2}}}^{2}}}
\]
\\
Repeating the procedure once again for the triple sum:
\[
\sum _{\ell_{{3}}=3}^{\infty} \left( \sum _{\ell_{{2}}=2}^{\ell_{{3}}-1} \left( \sum _
{\ell_{{1}}=1}^{\ell_{{2}-1}} \frac{1}{\ell_1^2 \cdot\ell_2^2 \cdot\ell_3^2}  \right)  \right) 
\]
with the corresponding MAPLE code
\begin{verbatim}
sum(sum(sum(1/(l[1]^2*l[2]^2*l[3]^2),
l[1]=1..l[2]-1),l[2]=2..l[3]-1),l[3]=3..infinity);
\end{verbatim}
Once more, we do not obtain the desired result $\frac{\pi^6}{7!}$, but a more complicated expression
\[
\sum _{l_{{3}}=3}^{\infty } \left( -\frac{1}{6}\,{\frac {{\pi }^{2}\Psi
 \left( 1,l_{{3}} \right) }{{l_{{3}}}^{2}}}+\frac{1}{6}\,{\frac {{\pi }^{2}
 \left( -1+1/6\,{\pi }^{2} \right) }{{l_{{3}}}^{2}}}+\sum _{l_{{2}}=2
}^{l_{{3}}-1}-{\frac {\Psi \left( 1,l_{{2}} \right) }{{l_{{2}}}^{2}{l_
{{3}}}^{2}}} \right)
\]
\\
With exception of  $\frac{\pi^2}{6}$ MAPLE does not generate rational multiples of $\pi$ powers.
$ \\ $
$ \\ $
Next we test how MAPLE is handling the infinite sums computing with floating-point numbers. In MAPLE the
\begin{verbatim}
Digits
\end{verbatim} environment variable controls the number of digits that Maple uses when making calculations with floating-point numbers.  Setting
\begin{verbatim}
Digits:=50;
\end{verbatim}
MAPLE computes with accuracy of 50 Digits. The command  
\begin{verbatim}evalf() 
\end{verbatim}
 numerically evaluates the expression. 
$ \\ $
$ \\ $
First computing power $\frac{\pi^4}{5!}$ and its corresponding double sum
\begin{verbatim}
Digits:=50; #number of digits that Maple uses
evalf(Pi^4/5!);
evalf(sum(sum(1/(l[1]^2*l[2]^2),
l[1]=1..l[2]-1),l[2]=2..infinity));
\end{verbatim}
$ \\ $
MAPLE presents
$ \\ $
\centerline{$Digits:=50$ } 
\centerline{$0.81174242528335364363700277240587592708106321393904$}
\centerline{$0.81174242528335364363700277240587592708106321393905$}
$ \\ $
The two values are very exact. They differ only in the last Digit.
$ \\ $
$ \\ $
Computing power $\frac{\pi^6}{7!}$ and comparing it with the infinite triple sum
\begin{verbatim}
evalf(Pi^6/7!);
evalf(sum(sum(sum(1/(l[1]^2*l[2]^2*l[3]^2),
l[1]=1..l[2]-1),l[2]=2..l[3]-1),l[3]=3..infinity));
\end{verbatim}
gives\\
\centerline{$0.19075182412208421369647211183579759898159077938116$}
\centerline{$0.19075182412208421369647211183579759898159077938116$}
$ \\ $
The two value are the same. There is no difference.
$ \\ $
$ \\ $
Once more computing power $\frac{\pi^8}{9!}$  and its corresponding quadrupel sum
\begin{verbatim}
evalf(Pi^8/9!);
evalf(sum(sum(sum(sum(1/(l[1]^2*l[2]^2*l[3]^2*l[4]^2),
l[1]=1..l[2]-1),l[2]=2..l[3]-1),l[3]=3..l[4]-1),l[4]=4..infinity));
\end{verbatim}
yields\\
\centerline{$0.026147847817654800504653261419496157949452103923173$}
\centerline{$0.026147847817654800504653261419496157949452103923173$}
$ \\ $
Both values are equal.
It's astonishing that MAPLE can compute infinite multi sums with such high accuracy.
\section{Appendix}
We present the proof of the theorem (section 2) that states, the product $\prod_{k=1}^{M}(1+x_{k}t)$ can be represented by a polynomial, whose coefficients consists of nested sums:
\begin{align*}
\prod_{k=1}^{M}(1+x_{k}t)=1
+&\left[\sum_{\ell_1=1}^{M}x{_{\ell_1}}\right] t
\\
+&\left[\sum _{\ell_{{2}}=2}^{M} \left( \sum _{\ell_{{1}}=1}^{\ell_{{2}}-1}x_{{\ell_{{1}}}
} x_{{\ell_{{2}}}} \right)\right] t^2
\\
+&\left[\sum _{\ell_{{3}}=3}^{M} \left( \sum _{\ell_{{2}}=2}^{\ell_{{3}}-1} \left( 
\sum _{\ell_{{1}}=1}^{\ell_{{2}}-1}x_{{\ell_{{1}}}} x_{{\ell_{{2}}}} x_{{\ell_{{3}}}} \right)  \right)\right] t^3
\\
&\qquad\qquad\qquad\qquad{\vdots}
\\
+&\left[\sum\limits_{{\ell_M} = M}^M {\left( {\sum\limits_{{\ell_{M - 1}} = M - 1}^{{\ell_M} - 1}  \cdots \left( {\sum\limits_{{\ell_2} = 2}^{{\ell_3} - 1} {\left( {\sum\limits_{{\ell_1} = 1}^{{\ell_{2} - 1}} x_{\ell_1} x_{\ell_2} \cdots x_{\ell_{M-1}} x_{\ell_M} } \right)} } \right) \cdots } \right)}\right]t^M
\end{align*}
$ \\ $
The proof of the theorem is accomplished by an induction argument on $M$.
$ \\ $
$ \\ $
In order to reduce paperwork, we replace the nested sums by $\sigma$-notation
\[
\sum_{\ell_1=1}^{M}x{_{\ell_1}}
=\sigma_{M,1}
\]
\[
\sum _{\ell_{2}=2}^M 
\left(
\sum _{\ell_1=1}^{\ell_2-1}
x_{\ell_1} x_{\ell_2}
\right)
=\sigma_{M,2}
\]
\[
\sum _{\ell_{3}=3}^M
\left(\sum _{\ell_2=2}^{\ell_3-1}
\left(\sum _{\ell_1=1}^{\ell_2-1}
x_{\ell_1}x_{\ell_2}x_{\ell_3}
\right)\right)
=\sigma_{M,3}
\]
\[\vdots\]
\[
\sum _{\ell_{k}=k}^M
\left(\sum _{\ell_{k-1}=k-1}^{\ell_k-1}
\ldots
\left(\sum _{\ell_2=2}^{\ell_3-1}
\left( \sum _{\ell_1=1}^{\ell_2-1}
x_{\ell_1}x_{\ell_2}\ldots x_{\ell_M-1}x_{\ell_M} \right)\right)
\ldots\right)
=\sigma_{M,k}
\]
$$
$$
\[\vdots\]
$$
$$
\[
\sum _{\ell_{M}=M}^M
\left(\sum _{\ell_{M-1}=M-1}^{\ell_M-1}
\ldots
\left(\sum _{\ell_2=2}^{\ell_3-1}
\left( \sum _{\ell_1=1}^{\ell_2-1}
x_{\ell_1}x_{\ell_2}\ldots x_{\ell_M-1}x_{\ell_M} \right)\right)
\ldots\right)
=\sigma_{M,M}
\]
\\
That gives the simple sum
\begin{align*}
\prod_{k=1}^{M}(1+x_{k}t)
=1
+\sigma_{M,1}\; t
+\sigma_{M,2}\; t^2
+\ldots
+\sigma_{M,k-1}\; t^{k-1}
+\sigma_{M,M}\; t^M
=1+\sum_{k=1}^{M}\sigma_{M,k}\; t^{k}
\end{align*}
\[\]
\textbf{Proof}
$ \\ $
\underline{Step 1:}
\\
Let $M = 1$. Then the left hand of the theorem is 
\[\prod_{k=1}^{1}(1+x_{1}t)=1+x_{1}t\]
and the right hand side of the theorem is 
\[
1+\sigma_{1,1}t=1+\left[\sum_{\ell_1=1}^{1}x{_{\ell_1}}\right]t=1+x_{1}t
\]
 which both equal 
$1 + x_1t$. So the theorem is true for $M = 1$.
$ \\ $
$ \\ $
\underline{Step 2:}
$ \\ $
We assume the theorem is true for $M = m \geq 1$, ie. that
\begin{align*}
\prod_{k=1}^{m}(1+x_{k}t)
=1+\sum_{k=1}^{m}\sigma_{m,k}\; t^{k}
\end{align*}
$ \\ $
\underline{Step 3:}
$ \\ $
We want to prove that the theorem is true for $m + 1$, ie. that
\begin{align*}
\prod_{k=1}^{m+1}(1+x_{k}t)
=1+\sum_{k=1}^{m+1}\sigma_{m+1,k}\; t^{k}
\end{align*}
$ \\ $
How do we get to step $3$ from step $2\,$? Let's have a look at the left-hand that we can write
\[
\prod_{k=1}^{m+1}\left(1+x_{k}t\right)
=\left(\prod_{k=1}^m\left(1+x_{k}t\right)\right)\left(1+x_{m+1}t\right)
\]
Replacing the product
\[
\prod_{k=1}^{m}\left(1+x_{k}t\right)
\]
by induction hypothesis
\[
1+\sum_{k=1}^{m}\sigma_{m,k}\; t^{k}
\]
and multiplying by
$\left(1+x_{m+1}\right)$ yields
\begin{align*}
\left(1+\sum_{k=1}^{m}\sigma_{m,k}\; t^{k}\right)
\left( 1+x_{m+1}t\right)
\end{align*}
$ \\ $
With some effort in expanding the brace, rearranging the terms, and gathering terms with equal powers, we achieve:
$ $
\begin{align*}
\left(1+\sum_{k=1}^{m}\sigma_{m,k}\; t^{k}\right)
\left( 1+x_{m+1}t\right)
&=1+\sum_{k=1}^{m}\sigma_{m,k}\; t^{k}
+x_{m+1}t+\left(\sum_{k=1}^{m}\sigma_{m,k}\; t^{k}\right)x_{m+1}t
\\
&=1+x_{m+1}t+\sum_{k=1}^{m}\sigma_{m,k}\; t^{k}
+\sum_{k=1}^{m}\sigma_{m,k}\; x_{m+1} \; t^{k+1}
\end{align*}
The next task requires a bit more effort. In order to simplifying it, we write
\[
A=x_{m+1}t+\sum_{k=1}^{m}\sigma_{m,k}\; t^{k}
\]
\[
B=\sum_{k=1}^{m}\sigma_{m,k}\; x_{m+1} \; t^{k+1};
\]
The right hand side simplifies to
\[
1+A+B
\]
First we split the term $A$
\begin{align*}
A=x_{m+1}t+\sum_{k=1}^{m}\sigma_{m,k}\; t^{k}
&=x_{m+1}t+\sum_{k=1}^{1}\sigma_{m,k}\; t^{k}
+\sum_{k=2}^{m}\sigma_{m,k}\; t^{k}
\end{align*}
reordering and simplifying the terms gives
\begin{align*}
A
=x_{m+1}t+\sigma_{m,1}\;t+\sum_{k=2}^{m}\sigma_{m,k}\;t^{k}
=\left(x_{m+1}+\sigma_{m,1}\right)t+\sum_{k=2}^{m}\sigma_{m,k}\;t^{k}\;.
\end{align*}
\\
Further holds
\[
\sigma_{m,1}
=\sum_{\ell_1=1}^{m}x_{\ell_1}\;.
\]
So, we finally obtain
\begin{align*}
A
=\left(x_{m+1}+\sum_{\ell_1=1}^{m}x_{\ell_1}\right)t+\sum_{k=2}^{m}\sigma_{m,k}\;t^{k}
=\sum_{\ell_1=1}^{m+1}x_{\ell_1}t+\sum_{k=2}^{m}\sigma_{m,k}\;t^{k}
\end{align*}
Let's focus on term
\[
B=\sum_{k=1}^{m}\sigma_{m,k}\; x_{m+1} \; t^{k+1}
\]
Splitting and shifting the index leads to 
\begin{align*}
B&=\sum _{k=1}^{m}\sigma_{{m,k}}\;x_{{m+1}}\;{t}^{k+1}
=\sum _{k=1}^{m-1}\sigma_{{m,k}}\;x_{{m+1}}\;{t}^{k+1}
+\sum _{k=m}^{m}\sigma_{{m,k}}\;x_{{m+1}}\;{t}^{k+1}\\
&=\sum _{k=1}^{m-1}\sigma_{{m,k}}\;x_{{m+1}}{t}^{k+1}
+\sigma_{m,m} \cdot x_{m+1}\;t^{m+1}\\
&=\sum _{k=2}^{m}\sigma_{{m,k-1}}\;x_{{m+1}}\;{t}^{k}
+\left(\prod_{k=1}^{m}x_{k}\right) x_{m+1}t^{m+1}\\
&=\sum _{k=2}^{m}\sigma_{{m,k-1}}\;x_{{m+1}}\;{t}^{k}
+\left(\prod_{k=1}^{m+1}\;x_{k}\right)t^{m+1}\\
&=\sum _{k=2}^{m}\sigma_{{m,k-1}}\;x_{{m+1}}\;{t}^{k}
+\sum _{k=m+1}^{m+1}\sigma_{{m,k}}\;{t}^{k+1}
\end{align*}
\\
So far, we have
\begin{align*}
&A=\sum_{\ell_1=1}^{m+1}x_{\ell_1}t+\sum_{k=2}^{m}\sigma_{m,k}\;t^{k}
\\
&B=\sum _{k=2}^{m}\sigma_{{m,k-1}}\;x_{{m+1}}\;{t}^{k}
+\sum _{k=m+1}^{m+1}\sigma_{{m,k}}\;{t}^{k+1}\;.
\end{align*}
\\Putting all the pieces together, we obtain
\begin{align*}
1+A+B
&=1+\sum _{k=1}^{1}\sigma_{{m+1,k}}\;{t}^{k}
+\sum _{k=2}^{m}\sigma_{{m,k}}\;{t}^{k}
+\sum _{k=2}^{m}\sigma_{{m,k-1}}x_{{m+1}}\;{t}^{k}
+\sum _{k=m+1}^{m+1}\sigma_{{m,k}}\;{t}^{k+1}\\
&=1+\sum _{k=1}^{1}\sigma_{{m+1,k}}\;{t}^{k}
+\sum _{k=2}^{m}\left(\sigma_{{m,k}}+\sigma_{{m,k-1}}x_{{m+1}}\right)\;{t}^{k}
+\sum _{k=m+1}^{m+1}\sigma_{{m,k}}\;{t}^{k+1}\end{align*}
\\
Now, applying Lemma 1 yields:
\[
\sum _{k=2}^{m}\left(\sigma_{{m,k}}+\sigma_{{m,k-1}}x_{{m+1}}\right)\;{t}^{k}=\sum_{k=2}^{m}\sigma_{m+1,k}\; t^{k}
\]
Plugging in
\[
1+\sum_{k=1}^{1}\sigma_{m+1,k}\; t^{k}
+\sum_{k=2}^{m}\sigma_{m+1,k}\; t^{k}
+\sum_{k=m+1}^{m+1}\sigma_{m+1,k}\; t^{k}
\]
Combining the three sums to a single sum, we obtain
\[
\left(1+\sum_{k=1}^{m}\sigma_{m,k}\; t^{k}\right)
\left( 1+x_{m+1}t\right)=1+\sum_{k=1}^{m+1}\sigma_{m+1,k}\; t^{k}
=\prod_{k=1}^{m+1}(1+x_{k}t)\]
This is the result as required. Hence the theorem holds true for all $ m\geq 1$ and the proof is completed.

\selectlanguage{english}

\newpage
\selectlanguage{ngerman}
\begin{center}
{\LARGE\bf
Reihen für geradzahlige Potenzen von Pi
}
\vskip 0.25cm
{\LARGE\bf
 durch Verallgemeinerung von Euler's
}
\vskip 0.25cm
{\LARGE\bf
Lösungsmethode für das Basler Problem
}
$$
$$
\large
{Alois Schiessl}
\vskip 0.25cm
{E-Mail: aloisschiessl@web.de}
\end{center}
$ \\ $
$ \\ $
\centerline{\textbf{Zusammenfassung}}
$ \\ $
In dieser Abhandlung stellen wir ein Verfahren vor, das die Berechnung von Reihen für geradzahlige Potenzen von $\pi$ ermöglicht. Die Grundidee ist eine Verallgemeinerung des Verfahrens von Euler, mit dem er 1735 das Basler Problem löste. Wir stellen elementar-symmetrische Polynome durch mehrfach verschachtelte Summen dar und leiten davon Reihen für geradzahlige Potenzen der Kreiszahl $\pi$ ab, wie zum Beispiel:
\[
\frac{\pi^2}{3!}=\sum _{\ell_{{1}}=1}^{\infty} \frac{1}{\ell_1^2}
\;;\quad
\frac{\pi^4}{5!} 
=\sum _{\ell_{{2}}=2}^{\infty} \left( \sum _{\ell_{{1}}=1}^{\ell_{{2}}-1} 
\frac{1}{\ell_1^2 \cdot\ell_2^2}\right) 
\;;\quad
\frac{\pi^6}{7!}=\sum _{\ell_{{3}}=3}^{\infty} \left( \sum _{\ell_{{2}}=2}^{\ell_{{3}}-1} \left( \sum _
{\ell_{{1}}=1}^{\ell_{{2}-1}} \frac{1}{\ell_1^2 \cdot\ell_2^2 \cdot\ell_3^2}  \right)  \right) 
\quad\cdots
\]
Die meisten der angegebenen Reihen scheinen nicht so bekannt zu sein.

\centerline{}
\centerline{\it Veröffentlicht anlässlich des Tages der Kreiszahl $\pi$}
\centerline{\it 3 - 14 - 2024}

\section*{1\quad Einführung}
Das Basler Problem war ein sehr berühmtes mathematisches Problem Mitte des 17. Jahrhunderts. Es wurde erstmals 1650 von Pietro Mengoli formuliert. Das Basler Problem besteht in der Aufgabe, für die Summe der Kehrwert-Quadrate einen geschlossenen Ausdruck zu finden. Es dauerte fast 100 Jahre bis Euler \cite{E1}, \cite{E2} im Jahre 1735 ein Ergebnis vorlegen konnte: 
\[
\sum_{k=1}^\infty \frac{1}{k^2}=\frac{\pi^2}{6}\;.
\]
\\
Euler's Idee beruht auf dem Zusammenhang von Polynomen und ihren Nullstellen. Ein Polynom vom Grad $n$ und führendem Koeffizienten 1 kann als Produkt $\left(x-x_1\right)\ldots\left(x-x_n\right)$ angegeben werden, wobei $x_1\ldots,x_n$ die Nullstellen des Polynoms sind. Euler nahm an, dass diese Eigenschaft auch auf Polynome mit unendlich vielen Nullstellen und somit unendlich vielen Faktoren zutrifft. Um Eulers Gedankengang nachverfolgen zu können, erinnern wir uns an die Potenzreihendarstellung von $\sin(x)$
$ $
\[ \sin(x)=x-\frac{x^3}{3!}+\frac{x^5}{5!}-+\cdots\]
Unter Verwendung dieser Potenzreihe definierte Euler die Funktion
\[
{P}\left(x\right)=\left\{
\begin{aligned}
\frac{\sin \left( x \right)}{x}\,, \quad &x \ne 0 \\ 
1\,,\;\quad\quad &x = 0\\ 
\end{aligned} 
\right.
\]
Da die Nullstellen von $P\left(x\right)$ durch $\{\pm\pi;\pm2\pi;\pm3\pi;\ldots\}$ gegeben sind, folgerte er weiter, dass $P\left(x\right)$ durch ein unendliches Produkt beschrieben werden könne: 
\[
\frac{\sin(x)}{x}=
\left(1-\frac{x^2}{\pi^2}\right)\left(1-\frac{x^2}{2^2\pi^2}\right)\left(1-\frac{x^2}{3^2\pi^2}\right)\ldots\]
Folge dessen müsste die folgende Darstellung Gültigkeit haben
\[
1-\frac{x^2}{3!}+\frac{x^4}{5!}-+\cdots=
\left(1-\frac{x^2}{\pi^2}\right)\left(1-\frac{x^2}{2^2\pi^2}\right)\left(1-\frac{x^2}{3^2\pi^2}\right)\ldots
\]
Da die Terme auf der linken und rechten Seite dieselben Nullstellen haben und überdies den gleichen Funktionswert  für $x=0$ annehmen, so folgerte Euler, handelt sich um zwei verschiedene Darstellungen ein und derselben Funktion. Leider gilt diese Schlussfolgerung nicht so ohne weiteres. Das ahnte Euler und viele seiner Zeitgenossen. Die komplexe Analysis war zu jener Zeit noch nicht weit genug entwickelt. Nichts desto trotz löste Euler die Klammern auf der rechten Seite auf, fasste die $x^2$ Terme zusammen und erhielt schließlich sein berühmtes Resultat
\[
\frac{\pi^2}{6}=\sum_{k=1}^\infty \frac{1}{k^2}\;.
\]
Erst hundert Jahre später bewies Karl Weierstrass mit Methoden der komplexen Analysis, dass Euler's Vorgehensweise richtig war.
$ \\ $
$ \\ $
In dieser Abhandlung verallgemeinern wir Eulers Methode und übertragen sie auf höhere Potenzen von $x$. Ausgangspunkt sind elementar-symmetrische Polynome, die wir durch mehrfache ineinander geschachtelte Summen ersetzen. Unter Einbeziehung der reellen $sinc(x)$  Funktion erhalten wir daraus dann die gesuchten Reihendarstellungen für geradzahlige Potenzen von $\pi$.

\section*{2\quad Ineinander geschachtelte Summen}
Zu $M\in\mathbb N$ betrachten wir die $M$ Unbestimmten $x_{1},\ldots ,x_{M}$ und bilden das Produkt
\begin{align*}
\left(1+x_1 t\right)\cdots\left(1+x_M t\right)=\prod_{k=1}^{M}(1+x_{k}t)
\end{align*}
Wir lösen die Klammern auf und fassen gleiche Potenzen zusammen, so dass wir schließlich die folgende Darstellung erhalten:
\begin{align*}
\prod_{k=1}^{M}(1+x_{k}t)
=1+\sigma_{M,1}t+\sigma_{M,2}t^2+\sigma_{M,3}t^3+\ldots+\sigma_{M,M}t^M=1+\sum_{k=1}^{M}\sigma_{M,k}\cdot t^{k}
\end{align*}
Hierbei bestehen die Koeffizienten aus den elementar-symmetrischen Polynomen
\begin{align*}
\sigma_{M,k}&=\begin{cases}
{\sum\limits_{1 \le {\ell_1} <\cdots < {\ell_k} \le M} {x_{{\ell_1}}\cdots x_{{\ell_k}}} } & (1\leq k\leq M) \\
\\
\qquad 0 & (M<k)
\end{cases}
\end{align*}
Wir schreiben sie einmal einzeln hin. Sei $k\in\mathbb \{1,2,3,...,M\}$, so ergeben sich der Reihe nach:
\[
\sigma_{M,1}=\sum\limits_{1 \le {\ell_1} \le M} {x{}_{{\ell_1}}}
\]
\[
\sigma_{M,2}=\sum\limits_{1 \le {\ell_1} < {\ell_2} \le M} 
x_{\ell_1}\cdot x_{\ell_2}
\]
\[
\sigma_{M,3}=\sum\limits_{1 \le {\ell_1} < {\ell_2} < {\ell_3} \le M} x_{\ell_1}\cdot x_{\ell_2}\cdot x_{\ell_3}
\]
\centerline{\vdots}
\[
\sigma_{M,k}={\sum\limits_{1 \le {\ell_1} <\cdots < {\ell_k} \le N} {x_{{\ell_1}}\cdots x_{{\ell_k}}} }
\]
\centerline{\vdots}
\[
\sigma_{M,M}={\sum\limits_{1 \le {\ell_1} <\cdots < {\ell_M} \le M} {x_{{\ell_1}}\cdots x_{{\ell_M}}} }=\prod_{\ell_k=1}^{M}x_{\ell_k}
\]
Wir geben ein Lemma \cite{EF} über Eigenschaften von elementar-symmetrischen Polynomen an, das sich später als sehr nützlich erweisen wird.
$ \\ $
$ \\ $
\textbf{Lemma 1}
\begin{align*}
\sigma_{M,k}&=0, \;M < k
\\
\sigma_{M,0}&=1, \;\; 0 \;\leq M
\\
\sigma_{M+1,k}&=\sigma_{M,k}+x_{M+1}\; \sigma_{M,k-1},\;\;0 < k \le M+1
\end{align*}
Die elementar-symmetrischen Polynome in der allgemeinen Darstellung
\[
{\sum\limits_{1 \le {\ell_1} <\cdots < {\ell_k} \le M} {x_{{\ell_1}}\cdots x_{{\ell_k}}} }
\]
sind mit ihren ineinander geschachtelten Muli-Indizes für das praktische Rechnen nicht so ohne weiteres verwendbar. Es wäre eine direktere Darstellung wünschenswert. Das sieht auf den ersten Blick kompliziert aus, ist aber im Prinzip ganz einfach. Die Idee hierzu ist die Umwandlung in mehrfach geschachtelte Summen. Man muss hierzu nur die Multi-Ungleichung  
\[
1 \le {\ell_1} < {\ell_2} <\cdots < {\ell_{k-1}} < {\ell_k} \le M
\]
in $k$ einzelne Ungleichungen auflösen, die einzelnen  Start- und Ende-Indizes feststellen, und sie dann durch $k$ ineinander geschachtelte Summen darstellen.

$ \\ $
Wir fangen mit $k=1$ an. Aus der Ungleichung $1 \le {\ell_1} \le M$ erhalten wir als Start-Index $1 \le\ell_1$  und als Ende-Index $\ell_1 \le M$. Das ergibt die einfache Summe
\[
\sigma_{M,1}
=\sum_{\ell_1=1}^{M}x_{\ell_1}
=x_1+\ldots+x_M
\]
\\
Nun sei $k=2$. Wir zerlegen $1 \le\ell_1 <\ell_2 \le M$ in zwei Ungleichungen
\begin{align*}
1 \le\ell_1 &\;\;\text{und}\;\;\ell_1 \le\ell_2-1\\
2 \le\ell_2 &\;\;\text{und}\;\;\ell_2 \le M
\end{align*}
Das ergibt dann eine Doppelsumme
\[
\sigma_{M,2}
=\underbrace {\sum _{\ell_{{2}}=2}^{M} \left( \sum _{\ell_{{1}}=1}^{\ell_{{2}}-1}x_{{\ell_{{1}}}
}\cdot x_{{\ell_{{2}}}} \right)}_{2\;-fach \;Summe}\\
=
x_1x_2+\ldots+x_{M-1}x_M
\]
\\
Im Fall $k=3$ bilden wir aus der Dreifach-Ungleichung  $1 \le\ell_1 <\ell_2 <\ell_3 \le M$ drei einzelne Ungleichungen
\begin{align*}
1 \le\ell_1 &\;\;\text{und}\;\;\ell_1 \le\ell_2-1\\
2 \le\ell_2 &\;\;\text{und}\;\;\ell_2 \le\ell_3-1\\
3 \le\ell_3 &\;\;\text{und}\;\;\ell_3 \le M
\end{align*}
und erhalten drei ineinander geschachtelte Summen
\[
\sigma_{M,3}
=\underbrace {\sum _{\ell_{{3}}=3}^{M} \left( \sum _{\ell_{{2}}=2}^{\ell_{{3}}-1} \left( 
\sum _{\ell_{{1}}=1}^{\ell_{{2}}-1} x_{{\ell_{{1}}}} \cdot x_{{\ell_{{2}}}}\cdot x_{{\ell_{{3}}}} \right)  \right)}_{3\;-fach\;Summe}
\]
\[
=x_1\cdot x_2\cdot x_3+\cdots+x_{M-2}\cdot x_{M-1}\cdot x_M
\]
\\
In derselben Weise ergeben sich für $k=4$
aus der Vierfach-Ungleichung  $1 \le\ell_1 <\ell_2 <\ell_3 <\ell_4 \le M$ vier einzelne Ungleichungen
\begin{align*}
1 \le\ell_1 &\;\;\text{und}\;\;\ell_1 \le\ell_2-1\\
2 \le\ell_2 &\;\;\text{und}\;\;\ell_2 \le\ell_3-1\\
3 \le\ell_3 &\;\;\text{und}\;\;\ell_3 \le\ell_4-1\\
4 \le\ell_4 &\;\;\text{und}\;\;\ell_4 \le M
\end{align*}
Das ergibt dann vier ineinander geschachtelte Summen
\[\sigma_{M,4}
=\underbrace{\sum _{\ell_{{4}}=4}^{M} \left(\sum _{\ell_{{3}}=3}^{\ell_4-1} \left( \sum _{\ell_{{2}}=2}^{\ell_{{3}}-1} \left( 
\sum _{\ell_{{1}}=1}^{\ell_{{2}}-1} x_{{\ell_{{1}}}} \cdot x_{{\ell_{{2}}}}\cdot x_{{\ell_3}}\cdot x_{{\ell_4}} \right)  \right)\right)}_{4\;-fach\;Summe}
\]
\[
=x_1\cdot x_2\cdot x_3\cdot x_4+\cdots+x_{M-3}\cdot x_{M-2}\cdot x_{M-1}\cdot x_M
\]
\\
Kommen wir zum allgemeinen Fall $k\le M$. Die $k$-fache Ungleichung 
\[1 \le {\ell_1} <\ell_2 \cdots <l_{k-1}  < {\ell_k} \le M\]
ergibt $k$ Ungleichungen 
\begin{align*}
1 \le\ell_1 &\;\;\text{und}\;\;\ell_1 \le\ell_2-1\\
2 \le\ell_2 &\;\;\text{und}\;\;\ell_2 \le\ell_3-1\\
&\quad\,\vdots\\
k-1 \le\ell_{k-1} &\;\;\text{und}\;\;\ell_{k-1} \le\ell_k-1\\
k \le\ell_k &\;\;\text{and}\;\;\ell_k \le M
\end{align*}
Hieraus erhalten wir $k$ ineinander geschachtelte Summen:
\[
\sigma_{M,k}
=\underbrace {\sum\limits_{{\ell_k} = k}^M {\left( {\sum\limits_{{\ell_{k - 1}} = k - 1}^{{\ell_k} - 1}  \cdots  \left( {\sum\limits_{{\ell_2} = 2}^{{\ell_3} - 1} {\left( {\sum\limits_{{\ell_1} = 1}^{{\ell_{2} - 1}} x_{\ell_1}\cdot x_{\ell_2} \cdots x_{\ell_{k-1}} \cdot x_{\ell_k} } \right)} } \right) \cdots } \right)} }_{k\;-fache\;Summe}
\]
\[={x_1\cdots x_k}+\cdots+x_{M-k+1}\cdots x_M 
\]\\
Anmerkung: Die Anzahl der Terme von $\sigma_{M,k}$ beträgt $\#\sigma_{M,k}=\binom{M}{k}.$
$ \\ $
$ \\ $
Bis jetzt haben wir folgende Reihen berechnet
\[
\prod_{k=1}^{1}(1+x_{k}t)=
1
+\left[\sum_{\ell_1=1}^{1}x{_{\ell_1}}\right] t
\]

\[
\prod_{k=1}^{2}(1+x_{k}t)=
1
+\left[\sum_{\ell_1=1}^{2}x{_{\ell_1}}\right] t
+\left[\sum _{\ell_{{2}}=2}^{2} \left( \sum _{\ell_{{1}}=1}^{\ell_{{2}}-1}x_{{\ell_{{1}}}
} x_{{\ell_{{2}}}} \right)\right] t^2
\]

\[
\prod_{k=1}^{3}(1+x_{k}t)=
1
+\left[\sum_{\ell_1=1}^{3}x{_{\ell_1}}\right] t
+\left[\sum _{\ell_{{2}}=2}^{3} \left( \sum _{\ell_{{1}}=1}^{\ell_{{2}}-1}x_{{\ell_{{1}}}
} x_{{\ell_{{2}}}} \right)\right]t^2
+\left[\sum _{\ell_{{3}}=3}^{3} \left( \sum _{\ell_{{2}}=2}^{\ell_{{3}}-1} \left( 
\sum _{\ell_{{1}}=1}^{\ell_{{2}}-1}x_{{\ell_{{1}}}} x_{{\ell_{{2}}}} x_{{\ell_{{3}}}} \right)  \right)\right] t^3
\]
$ \\ $
Für $M\in\mathbb N$ formulieren wir folgenden Lehrsatz:
$ \\ $
$ \\ $
\textbf{Theorem 1}
\textit{Das Produkt}
\begin{align*}
\prod_{k=1}^{M}(1+x_{k}t)=\left(1+x_1 t\right)\cdots\left(1+x_M t\right)
\end{align*}
\textit{kann durch ein Polynom dargestellt werden, dessen Koeffizienten aus den folgenden mehrfach ineinander geschachtelten Summen bestehen:}
\begin{align*}
\prod_{k=1}^{M}(1+x_{k}t)=1
+&\left[\sum_{\ell_1=1}^{M}x{_{\ell_1}}\right] t
\\
+&\left[\sum _{\ell_{{2}}=2}^{M} \left( \sum _{\ell_{{1}}=1}^{\ell_{{2}}-1}x_{{\ell_{{1}}}
} x_{{\ell_{{2}}}} \right)\right] t^2
\\
+&\left[\sum _{\ell_{{3}}=3}^{M} \left( \sum _{\ell_{{2}}=2}^{\ell_{{3}}-1} \left( 
\sum _{\ell_{{1}}=1}^{\ell_{{2}}-1}x_{{\ell_{{1}}}} x_{{\ell_{{2}}}} x_{{\ell_{{3}}}} \right)  \right)\right] t^3
\\
+&\left[{\sum _{\ell_{{4}}=4}^{M} \left(\sum _{\ell_{{3}}=3}^{\ell_4-1} \left( \sum _{\ell_{{2}}=2}^{\ell_{{3}}-1} \left( 
\sum _{\ell_{{1}}=1}^{\ell_{{2}}-1} x_{{\ell_{{1}}}}  x_{{\ell_{{2}}}} x_{{\ell_3}} x_{{\ell_4}} \right)  \right)\right)}\right] t^4
\\
&\qquad\qquad\qquad\qquad{\vdots}
\\
+&\left[\sum\limits_{{\ell_M} = M}^M {\left( {\sum\limits_{{\ell_{M - 1}} = M - 1}^{{\ell_M} - 1}  \cdots \left( {\sum\limits_{{\ell_2} = 2}^{{\ell_3} - 1} {\left( {\sum\limits_{{\ell_1} = 1}^{{\ell_{2} - 1}} x_{\ell_1} x_{\ell_2} \cdots x_{\ell_{M-1}} x_{\ell_M} } \right)} } \right) \cdots } \right)}\right]t^M
\end{align*}
Der Beweis erfolgt mit vollständiger Induktion nach $M$. Der Beweis ist nicht schwierig, aber mit etwas Aufwand verbunden, so dass wir ihn am Ende der Abhandlung bringen.

\section*{3\quad Von der Polynomdarstellung zur sinc function}
In der Darstellung
\begin{align*}
\prod_{k=1}^{M}(1+x_{k}t)=1
+&\left[\sum_{\ell_1=1}^{M}x{_{\ell_1}}\right] t
\\
+&\left[\sum _{\ell_{{2}}=2}^{M} \left( \sum _{\ell_{{1}}=1}^{\ell_{{2}}-1}x_{{\ell_{{1}}}
} x_{{\ell_{{2}}}} \right)\right] t^2
\\
+&\left[\sum _{\ell_{{3}}=3}^{M} \left( \sum _{\ell_{{2}}=2}^{\ell_{{3}}-1} \left( 
\sum _{\ell_{{1}}=1}^{\ell_{{2}}-1}x_{{\ell_{{1}}}} x_{{\ell_{{2}}}} x_{{\ell_{{3}}}} \right)  \right)\right] t^3
\\
+&\left[{\sum _{\ell_{{4}}=4}^{M} \left(\sum _{\ell_{{3}}=3}^{\ell_4-1} \left( \sum _{\ell_{{2}}=2}^{\ell_{{3}}-1} \left( 
\sum _{\ell_{{1}}=1}^{\ell_{{2}}-1} x_{{\ell_{{1}}}}  x_{{\ell_{{2}}}} x_{{\ell_3}} x_{{\ell_4}} \right)  \right)\right)}\right] t^4
\\
&\qquad\qquad\qquad\qquad{\vdots}
\\
+&\left[\sum\limits_{{\ell_M} = M}^M {\left( {\sum\limits_{{\ell_{M - 1}} = M - 1}^{{\ell_M} - 1}  \cdots \left( {\sum\limits_{{\ell_2} = 2}^{{\ell_3} - 1} {\left( {\sum\limits_{{\ell_1} = 1}^{{\ell_{2} - 1}} x_{\ell_1} x_{\ell_2} \cdots x_{\ell_{M-1}} x_{\ell_M} } \right)} } \right) \cdots } \right)}\right]t^M
\end{align*}
 substituieren wir auf der linken Seite $t=-x^2$ und $x_k=\frac{1}{k^2}$. Damit die Gleichung richtig bleibt müssen wir auf der rechten Seite ebenfalls $t=-x^2$ and $x_{\ell_k}=\frac{1}{\ell_k^2}$ substituieren. Das ergibt dann
\begin{align*}
\prod_{k=1}^{M}(1-\frac{x^2}{k^2})=1-&\left[\sum_{\ell_1=1}^{N}\frac{1}{\ell_1^2} \right] x^2\\
+&\left[\sum _{\ell_{{2}}=2}^{N} \left( \sum _{\ell_{{1}}=1}^{\ell_{{2}}-1}\frac{1}{\ell_1^2\cdot \ell_2^2} \right)\right] x^4\\
-&\left[\sum _{\ell_{{3}}=3}^{N} \left( \sum _{\ell_{{2}}=2}^{\ell_{{3}}-1} \left( 
\sum _{\ell_{{1}}=1}^{\ell_{{2}}-1}
\frac{1}{\ell_1^2\cdot \ell_2^2\cdot \ell_3^2} \right)  \right)\right] x^6\\
+&\left[\sum _{\ell_{{4}}=4}^{N} \left( \sum _{\ell_{{3}}=3}^{\ell_{{4}}-1}
 \left( \sum _{\ell_{{2}}=2}^{\ell_{{3}}-1} \left( \sum _{\ell_{{1}}=1}^{\ell_{{2}
}-1}
\frac{1}{\ell_1^2\cdot \ell_2^2\cdot \ell_3^2\cdot \ell_4^2} \right) 
 \right)  \right)\right] x^8
\\
&\qquad\qquad\qquad\qquad{\vdots}
\\
+\left(-1\right)^M
&\left[\sum\limits_{{\ell_M} = M}^M {\left( {\sum\limits_{{\ell_{M - 1}} = M - 1}^{{\ell_M} - 1}  \cdots  \left( {\sum\limits_{{\ell_2} = 2}^{{\ell_3} - 1} {\left( {\sum\limits_{{\ell_1} = 1}^{{\ell_{2 - 1}}} \frac{1}{\ell_1^2\cdot \ell_2^2 \cdots\ell_{M-1}^2\cdot \ell_M^2} } \right)} } \right) \cdots } \right)}\right]x^{2M}
\end{align*}
\[\]
Als nächstes überlegen wir uns, was geschieht, wenn wir  $M \to \infty $ streben lassen.
$ \\ $
Auf der linken Seite erhalten wir ein unendliches Produkt
\[
\underset{M\to \infty }{\mathop{\lim }}\,\prod_{k=1}^{M}(1-\frac{x^2}{k^2})
=\prod_{k=1}^{\infty}(1-\frac{x^2}{k^2})\,,
\]
das wir als die Produktdarstellung der $sinc$ Funktion identifizieren. Die $sinc$ Funktion ist definiert als: 
\[
\text{sinc}\left( \pi x \right)=\left\{ \begin{aligned}
\frac{\sin \left( \pi x \right)}{\pi x}, \qquad &x \ne 0 \\ 
1\;\,,\quad\qquad &x = 0
\end{aligned} 
\right.
\]
$ \\ $
Zu ihr gehört die Potenzreihe, die sich leicht aus der Potenzreihe $sinus$ Funktion ableiten lässt:
\[
\frac{\sin\left(\pi x\right)}{\pi x}=1-\frac{\pi^2}{3!}x^2+\frac{\pi^4}{5!}x^4-\frac{\pi^6}{7!}x^6+-\cdots
\]
$ \\ $
Auf der rechten Seite erhalten wir für $M \to \infty $ eine Potenzreihe, deren Koeffizienten selbst aus unendlichen Summen bestehen.
\begin{align*}
1-&\left[\sum_{\ell_1=1}^{\infty}\frac{1}{\ell_1^2} \right] x^2\\
+&\left[\sum _{\ell_{{2}}=2}^{\infty} \left( \sum _{\ell_{{1}}=1}^{\ell_{{2}}-1}\frac{1}{\ell_1^2\cdot \ell_2^2} \right)\right] x^4\\
-&\left[\sum _{\ell_{{3}}=3}^{\infty} \left( \sum _{\ell_{{2}}=2}^{\ell_{{3}}-1} \left( 
\sum _{\ell_{{1}}=1}^{\ell_{{2}}-1}
\frac{1}{\ell_1^2\cdot \ell_2^2\cdot \ell_3^2} \right)  \right)\right] x^6\\
+&\left[\sum _{\ell_{{4}}=4}^{\infty} \left( \sum _{\ell_{{3}}=3}^{\ell_{{4}}-1}
 \left( \sum _{\ell_{{2}}=2}^{\ell_{{3}}-1} \left( \sum _{\ell_{{1}}=1}^{\ell_{{2}
}-1}
\frac{1}{\ell_1^2\cdot \ell_2^2\cdot \ell_3^2\cdot \ell_4^2} \right) 
 \right)  \right)\right] x^8-+\ldots
\end{align*}

\section*{4 \quad Von der sinc Funktion zu den $\pi$-Potenzen}
Wir haben somit zwei Potenzreihendarstellungen für die $sinc$ Funktion. Einmal die klassische
\[
\frac{\sin\left(\pi x\right)}{\pi x}=1-\frac{\pi^2}{3!}x^2+\frac{\pi^4}{5!}x^4-\frac{\pi^6}{7!}x^6+-\cdots
\]
und zum anderen die Darstellung durch mehrfache Summen 
\begin{align*}
\frac{\sin\left(\pi x\right)}{\pi x}=1-&\left[\sum_{\ell_1=1}^{\infty}\frac{1}{\ell_1^2} \right] x^2\\
+&\left[\sum _{\ell_{{2}}=2}^{\infty} \left( \sum _{\ell_{{1}}=1}^{\ell_{{2}}-1}\frac{1}{\ell_1^2\cdot \ell_2^2} \right)\right] x^4\\
-&\left[\sum _{\ell_{{3}}=3}^{\infty} \left( \sum _{\ell_{{2}}=2}^{\ell_{{3}}-1} \left( 
\sum _{\ell_{{1}}=1}^{\ell_{{2}}-1}
\frac{1}{\ell_1^2\cdot \ell_2^2\cdot \ell_3^2} \right)  \right)\right] x^6\\
+&\left[\sum _{\ell_{{4}}=4}^{\infty} \left( \sum _{\ell_{{3}}=3}^{\ell_{{4}}-1}
 \left( \sum _{\ell_{{2}}=2}^{\ell_{{3}}-1} \left( \sum _{\ell_{{1}}=1}^{\ell_{{2}
}-1}
\frac{1}{\ell_1^2\cdot \ell_2^2\cdot \ell_3^2\cdot \ell_4^2} \right) 
 \right)  \right)\right] x^8-+\ldots
\end{align*}
Gemäß dem Eindeutigkeitssatz für Potenzreihen von Funktionen müssen beide identisch sein und somit ihre Koeffizienten übereinstimmen.
Es muss also gelten
\begin{align*}
1-&\frac{\pi^2}{3!}x^2+\frac{\pi^4}{5!}x^4-\frac{\pi^6}{7!}x^6+\frac{\pi^8}{9!}x^8-+\cdots
\\
\\
=1-&\left[\sum_{\ell_1=1}^{\infty}\frac{1}{\ell_1^2} \right] x^2\\
+&\left[\sum _{\ell_{{2}}=2}^{\infty} \left( \sum _{\ell_{{1}}=1}^{\ell_{{2}}-1}\frac{1}{\ell_1^2\cdot \ell_2^2} \right)\right] x^4\\
-&\left[\sum _{\ell_{{3}}=3}^{\infty} \left( \sum _{\ell_{{2}}=2}^{\ell_{{3}}-1} \left( 
\sum _{\ell_{{1}}=1}^{\ell_{{2}}-1}
\frac{1}{\ell_1^2\cdot \ell_2^2\cdot \ell_3^2} \right)  \right)\right] x^6\\
+&\left[\sum _{\ell_{{4}}=4}^{\infty} \left( \sum _{\ell_{{3}}=3}^{\ell_{{4}}-1}
 \left( \sum _{\ell_{{2}}=2}^{\ell_{{3}}-1} \left( \sum _{\ell_{{1}}=1}^{\ell_{{2}
}-1}
\frac{1}{\ell_1^2\cdot \ell_2^2\cdot \ell_3^2\cdot \ell_4^2} \right) 
 \right)  \right)\right] x^8-+\ldots
\end{align*}
$ \\ $
$ \\ $
Setzen wir die Koeffizienten von $x^2$ gleich, so erhalten wir
\[
\frac{\pi^2}{3!}=\sum _{\ell_{{1}}=1}^{\infty} \frac{1}{\ell_1^2} 
\]
Das ist das berühmte Ergebnis, das Euler 1735 publizierte. Es ist verwunderlich warum Euler nicht höhere Potenzen untersuchte. Vermutlich war er nur an Reihen der Form
\[
\sum_{k=1}^\infty \frac{1}{k^n}
\;\text{,}\;\; n=2,3,\dots
\]
interessiert. Wir tun es und fahren mit höheren Potenzen fort.
$ \\ $
$ \\ $
Setzen wir die Koeffizienten von $x^4$ gleich, so erhalten wir auf der linken Seite ein rationales Vielfaches von $\pi^4$ and auf der rechten Seite zwei ineinander geschachtelte Summen
\[
\frac{\pi^4}{5!} 
=\sum _{\ell_{{2}}=2}^{\infty} \left( \sum _{\ell_{{1}}=1}^{\ell_{{2}}-1} 
\frac{1}{\ell_1^2 \cdot\ell_2^2}\right) 
\]
Hierbei wollen wir den Grenzwert
\[
\sum _{\ell_{{2}}=2}^{\infty} \left( \sum _{\ell_{{1}}=1}^{\ell_{{2}}-1} 
\frac{1}{\ell_1^2 \cdot\ell_2^2}\right)
\]
$ \\ $
als den Grenzwert der Partialsummen verstehen
\[
\sum _{\ell_{{2}}=2}^{2} \left( \sum _{\ell_{{1}}=1}^{\ell_{{2}}-1}{\frac {1}{{\ell_{{1}}}^{2}{\ell_{{2}}}^{2}}} \right) 
=
\frac{1}{1^2\cdot 2^2}
\]
\[
\sum _{\ell_{{2}}=2}^{3} \left( \sum _{\ell_{{1}}=1}^{\ell_{{2}}-1}{\frac {1}{{\ell_{{1}}}^{2}{\ell_{{2}}}^{2}}} \right) 
=\frac{1}{1^2\cdot 2^2}+\frac{1}{1^2\cdot 3^2}+\frac{1}{2^2\cdot 3^2}
\]
\[
\sum _{\ell_{{2}}=2}^{4} \left( \sum _{\ell_{{1}}=1}^{\ell_{{2}}-1}{\frac {1}{{\ell_{{1}}}^{2}{\ell_{{2}}}^{2}}} \right) 
=
\frac{1}{1^2\cdot 2^2}+\frac{1}{1^2\cdot 3^2}+\frac{1}{1^2\cdot 4^2}+
\frac{1}{2^2\cdot 3^2}+\frac{1}{2^2\cdot 4^2}+
\frac{1}{3^2\cdot 4^2} 
\]
\[
\vdots
\]
$ $
Für die Koeffizienten von $x^6$ erhalten wir ebenfalls ein  rationales Vielfaches von  $\pi^6$, dargestellt durch eine dreifach verschachtelte Summe
\[
\frac{\pi^6}{7!}=\sum _{\ell_{{3}}=3}^{\infty} \left( \sum _{\ell_{{2}}=2}^{\ell_{{3}}-1} \left( \sum _
{\ell_{{1}}=1}^{\ell_{{2}-1}} \frac{1}{\ell_1^2 \cdot\ell_2^2 \cdot\ell_3^2}  \right)  \right) 
\]
Auch hier definieren wir den Grenzwert
\[
\sum _{\ell_{{3}}=3}^{\infty} \left( \sum _{\ell_{{2}}=2}^{\ell_{{3}}-1} \left( \sum _
{\ell_{{1}}=1}^{\ell_{{2}-1}} \frac{1}{\ell_1^2 \cdot\ell_2^2 \cdot\ell_3^2}  \right)  \right)
\]
als den Grenzwert der Partialsummen
\[
\sum _{\ell_{{3}}=3}^{3}
\left( \sum _{\ell_{{2}}=2}^{\ell_{{3}}-1}
\left( \sum _{\ell_{{1}}=1}^{\ell_{{2}}-1}{\frac {1}{{\ell_{{1}}}^{2}{\ell_{{2}}}^{2}{\ell_{{3}}}^{2}}} \right)  \right)
=
 \frac{1}{1^2 \cdot 2^2 \cdot 3^2}
\]

\begin{align*}
\sum _{\ell_{{3}}=3}^{4}
\left( \sum _{\ell_{{2}}=2}^{\ell_{{3}}-1}
\left( \sum _{\ell_{{1}}=1}^{\ell_{{2}}-1}{\frac {1}{{\ell_{{1}}}^{2}{\ell_{{2}}}^{2}{\ell_{{3}}}^{2}}} \right)  \right)
&=
 \frac{1}{1^2 \cdot 2^2 \cdot 3^2}
+\frac{1}{1^2 \cdot 2^2 \cdot 4^2}
+\frac{1}{1^2 \cdot 3^2 \cdot 4^2}
+\frac{1}{2^2 \cdot 3^2 \cdot 4^2}
\end{align*}

\begin{align*}
\sum _{\ell_{{3}}=3}^{5}
\left( \sum _{\ell_{{2}}=2}^{\ell_{{3}}-1}
\left( \sum _{\ell_{{1}}=1}^{\ell_{{2}}-1}{\frac {1}{{\ell_{{1}}}^{2}{\ell_{{2}}}^{2}{\ell_{{3}}}^{2}}} \right)  \right)
&=
 \frac{1}{1^2 \cdot 2^2 \cdot 3^2}
+\frac{1}{1^2 \cdot 2^2 \cdot 4^2}
+\frac{1}{1^2 \cdot 2^2 \cdot 5^2}
+\frac{1}{1^2 \cdot 3^2 \cdot 4^2}
+\frac{1}{1^2 \cdot 3^2 \cdot 5^2}\\
&+\frac{1}{1^2 \cdot 4^2 \cdot 5^2}
+\frac{1}{2^2 \cdot 3^2 \cdot 4^2}
+\frac{1}{2^2 \cdot 3^2 \cdot 5^2}
+\frac{1}{2^2 \cdot 4^2 \cdot 5^2}
+\frac{1}{3^2 \cdot 4^2 \cdot 5^2}
\end{align*}
\[
\vdots
\]
Für die weiteren Potenzen erhalten wir
\[
\frac{{{\pi ^8}}}{{9!}} =\sum\limits_{{\ell _4} = 4}^{\infty} {\left( {\sum\limits_{{\ell _3} = 3}^{{\ell _4} - 1} {\left( {\sum\limits_{{\ell _2} = 2}^{{\ell _3} - 1} {\left( {\sum\limits_{{\ell _1} = 1}^{{\ell _{2 - 1}}} {\frac{1}{{\ell _1^2 \cdot \ell _2^2 \cdot \ell _3^2 \cdot \ell _4^2}}} } \right)} } \right)} } \right)}
\]
\[
\frac{\pi^{10}}{11!}=\sum _{\ell_{{5}}=5}^{\infty}\left(\sum _{\ell_{{4}}=4}^{\ell_{{5}}-1}\left(\sum _{\ell_{{3}}=3}^{\ell_{{4}}-1} \left( \sum _{\ell_{{2}}=2}^{\ell_{{3}}-1} \left( \sum_
{\ell_{{1}}=1}^{\ell_{{2}-1}} \frac{1}{\ell_1^2 \cdot\ell_2^2 \cdot\ell_3^2 \cdot\ell_4^2 \cdot\ell_5^2}  \right)  \right) \right) \right)
\]
\[
\frac{\pi^{12}}{13!}=\sum _{\ell_{{6}}=6}^{\infty}\left(\sum _{\ell_{{5}}=5}^{\ell_{{6}}-1}
\left(\sum _{\ell_{{4}}=4}^{\ell_{{5}}-1}
\left(\sum _{\ell_{{3}}=3}^{\ell_{{4}}-1}
\left(\sum _{\ell_{{2}}=2}^{\ell_{{3}}-1}
\left( \sum_{\ell_{{1}}=1}^{\ell_{{2}-1}}
\frac{1}{\ell_1^2 \cdot\ell_2^2 \cdot\ell_3^2 \cdot\ell_4^2 \cdot\ell_5^2\cdot\ell_6^2}
\right)
\right)
\right)
\right)
\right)
\]
\[
\frac{\pi^{14}}{15!}=\sum_{\ell_{{7}}=7}^{\infty}\left(\sum _{\ell_{{6}}=6}^{\ell_{{7}}-1}
\left(\sum _{\ell_{{5}}=5}^{\ell_{{6}}-1}
\left(\sum _{\ell_{{4}}=4}^{\ell_{{5}}-1}
\left(\sum _{\ell_{{3}}=3}^{\ell_{{4}}-1}
\left(\sum _{\ell_{{2}}=2}^{\ell_{{3}}-1}
\left( \sum_{\ell_{{1}}=1}^{\ell_{{2}-1}}
\frac{1}{\ell_1^2 \cdot\ell_2^2 \cdot\ell_3^2 \cdot\ell_4^2 \cdot\ell_5^2\cdot\ell_6^2}
\right)
\right)
\right)
\right)
\right)
\right)
\]
\[
\vdots
\]
Im allgemeinen Falls $n\in\mathbb N$ ergibt sich
\[
\frac{{{\pi ^{2n}}}}{{\left( {2n + 1} \right)!}} = \sum\limits_{{\ell_n} = n}^{\infty} {\left( {\sum\limits_{{\ell_{n - 1}} = n - 1}^{{\ell_n} - 1}  \cdots  \left( {\sum\limits_{{\ell_2} = 2}^{{\ell_3} - 1} {\left( {\sum\limits_{{\ell_1} = 1}^{{\ell_{2 - 1}}} {\frac{1}{{\ell_1^2 \cdots\ell_n^2}}} } \right)} } \right) \cdots } \right)}
\]
Unter dem Grenzwert der Summe
\[
\sum\limits_{{\ell_n} = n}^{\infty} {\left( {\sum\limits_{{\ell_{n - 1}} = n - 1}^{{\ell_n} - 1}  \cdots  \left( {\sum\limits_{{\ell_2} = 2}^{{\ell_3} - 1} {\left( {\sum\limits_{{\ell_1} = 1}^{{\ell_{2 - 1}}} {\frac{1}{{\ell_1^2 \cdots\ell_n^2}}} } \right)} } \right) \cdots } \right)} 
\]
wollen wir auch hier den Grenzwert der Partialsummen zu verstehen
\[
\sum\limits_{{\ell_n} = n}^{n} {\left( {\sum\limits_{{\ell_{n - 1}} = n - 1}^{{\ell_n} - 1}  \cdots  \left( {\sum\limits_{{\ell_2} = 2}^{{\ell_3} - 1} {\left( {\sum\limits_{{\ell_1} = 1}^{{\ell_{2 - 1}}} {\frac{1}{{\ell_1^2 \cdots\ell_n^2}}} } \right)} } \right) \cdots } \right)} 
\]
\[
\sum\limits_{{\ell_n} = n}^{n+1} {\left( {\sum\limits_{{\ell_{n - 1}} = n - 1}^{{\ell_n} - 1}  \cdots  \left( {\sum\limits_{{\ell_2} = 2}^{{\ell_3} - 1} {\left( {\sum\limits_{{\ell_1} = 1}^{{\ell_{2 - 1}}} {\frac{1}{{\ell_1^2 \cdots\ell_n^2}}} } \right)} } \right) \cdots } \right)} 
\]
\[
\sum\limits_{{\ell_n} = n}^{n+2} {\left( {\sum\limits_{{\ell_{n - 1}} = n - 1}^{{\ell_n} - 1}  \cdots  \left( {\sum\limits_{{\ell_2} = 2}^{{\ell_3} - 1} {\left( {\sum\limits_{{\ell_1} = 1}^{{\ell_{2 - 1}}} {\frac{1}{{\ell_1^2 \cdots\ell_n^2}}} } \right)} } \right) \cdots } \right)} 
\]
\[
\vdots
\]
Natürlich müssen wir uns auch um die Konvergenz der Reihen kümmern. Zunächst stellen wir fest, dass es nur positive Summenden gibt. Somit empfiehlt sich das Majorantenkriterium das Mittel der Wahl. Wir müssen also eine Reihe finden deren Summenden alle positiv, stets größer oder höchsten gleich der untersuchten Reihe sind und die konvergiert.
$ \\ $
$ \\ $
Die Konvergenz der Reihe
\[
\sum _{\ell_{{1}}=1}^{\infty} \frac{1}{\ell_1^2} 
\]
ist wohl bekannt und wir haben nichts zu tun.
\\
Für höhere Potenzen verwenden wir die Vergleichsreihen \cite{KK}
\[
\sum _{\ell_{{2}}=1}^{\infty} \left( \sum _{\ell_{{1}}=1}^{\infty} 
\frac{1}{\ell_1^2 \cdot\ell_2^2}\right)=\left(\frac{\pi^2}{6}\right)^2;
\]
\[
\sum _{\ell_{{3}}=1}^{\infty} \left( \sum _{\ell_{{2}}=1}^{\infty} \left( \sum _
{\ell_{{1}}=1}^{\infty} \frac{1}{\ell_1^2 \cdot\ell_2^2 \cdot\ell_3^2}  \right)  \right)=\left(\frac{\pi^2}{6}\right)^3;
\]

\[
\vdots
\]

\[
\sum\limits_{{\ell_n} = 1}^{\infty} {\left( {\sum\limits_{{\ell_{n - 1}} =1}^{\infty}  \cdots  \left( {\sum\limits_{{\ell_2} = 1}^{\infty} {\left( {\sum\limits_{{\ell_1} = 1}^{\infty} {\frac{1}{{\ell_1^2 \cdots\ell_n^2}}} } \right)} } \right) \cdots } \right)}
=\left(\frac{\pi^2}{6}\right)^n\]
deren Grenzwert wohl bekannt ist und vergleichen sie mit den hergeleiteten Reihen. Das ergibt
\[
\sum _{\ell_{{2}}=2}^{\infty} \left( \sum _{\ell_{{1}}=1}^{\ell_{{2}}-1} 
\frac{1}{\ell_1^2 \cdot\ell_2^2}\right)
<\sum _{\ell_{{2}}=1}^{\infty} \left( \sum _{\ell_{{1}}=1}^{\infty} 
\frac{1}{\ell_1^2 \cdot\ell_2^2}\right)=\left(\frac{\pi^2}{6}\right)^2;
\]
\[
\sum _{\ell_{{3}}=3}^{\infty} \left( \sum _{\ell_{{2}}=2}^{\ell_{{3}}-1} \left( \sum _
{\ell_{{1}}=1}^{\ell_{{2}-1}} \frac{1}{\ell_1^2 \cdot\ell_2^2 \cdot\ell_3^2}  \right)  \right) 
<
\sum _{\ell_{{3}}=1}^{\infty} \left( \sum _{\ell_{{2}}=1}^{\infty} \left( \sum _
{\ell_{{1}}=1}^{\infty} \frac{1}{\ell_1^2 \cdot\ell_2^2 \cdot\ell_3^2}  \right)  \right)=\left(\frac{\pi^2}{6}\right)^3;
\]
\[
\vdots
\]
Im allgemeinen Fall für $n\in\mathbb N$ haben wir
\begin{align*}
&\sum\limits_{{\ell_n} = n}^{\infty} {\left( {\sum\limits_{{\ell_{n - 1}} = n - 1}^{{\ell_n} - 1}  \cdots  \left( {\sum\limits_{{\ell_2} = 2}^{{\ell_3} - 1} {\left( {\sum\limits_{{\ell_1} = 1}^{{\ell_{2 - 1}}} {\frac{1}{{\ell_1^2 \cdots\ell_n^2}}} } \right)} } \right) \cdots } \right)}\\
<
&\sum\limits_{{\ell_n} = 1}^{\infty} {\left( {\sum\limits_{{\ell_{n - 1}} =1}^{\infty}  \cdots  \left( {\sum\limits_{{\ell_2} = 1}^{\infty} {\left( {\sum\limits_{{\ell_1} = 1}^{\infty} {\frac{1}{{\ell_1^2 \cdots\ell_n^2}}} } \right)} } \right) \cdots } \right)}
=\left(\frac{\pi^2}{6}\right)^n
\end{align*}
Wir erkennen, dass jede $n$-fach ineinander geschachtelte Reihe durch $\left(\frac{\pi^2}{6}\right)^n $ beschränkt ist und hieraus folgt die Konvergenz der Reihen.
$ \\ $
$ \\ $
Anmerkung: Man kann die Konvergenz der Reihen auch damit begründen, dass bei den Summanden der Zähler stets 1 und im Nenner nur aufsteigende quadratische Terme der natürlichen Zahlen vorkommen. Die Reihe $\sum _{\ell_{{1}}=1}^{\infty} \frac{1}{\ell_1^2}$ kann dann als Majorante herangezogen werden.

\section*{5 \quad Berechnungen mit MAPLE}
Die Berechnungen wurden mit MAPLE Version 2024.0 vorgenommen.
\\
Falls die Software MAPLE zur Verfügung steht, können die Programmzeilen direkt per coppy \& paste in eine Maple Session eingefügt und ausgeführt werden.
$ \\ $
$ \\ $
Wir beginnen mit der Summe
\[
\sum _{\ell_{{1}}=1}^{\infty} \frac{1}{\ell_1^2} 
\]
Der MAPLE Befehl hierfür lautet

\begin{verbatim}
sum(1/l[1]^2,l[1]=1..infinity);
\end{verbatim}
Kopieren wir ihn per copy \& paste in eine MAPLE Session, so erhalten wir das erwartete Ergebnis
\[
\frac{1}{6}\,{\pi }^{2}
\]
Die Doppelsumme
\[
\sum _{\ell_{{2}}=2}^{\infty} \left( \sum _{\ell_{{1}}=1}^{\ell_{{2}}-1} 
\frac{1}{\ell_1^2 \cdot\ell_2^2}\right)
\]
schreibt sich in MAPLE so
\begin{verbatim}
sum(sum(1/(l[1]^2*l[2]^2),
l[1]=1..l[2]-1),l[2]=2..infinity);
\end{verbatim}
Ausgeführt gibt MAPLE nicht den erwarteten Wert $\frac{\pi^4}{5!}$ sondern einen Ausdruck mit der $\Psi$-Funktion aus:
\[
\sum _{l_{{2}}=2}^{\infty }-{\frac {\Psi \left( 1,l_{{2}} \right) }{{l
_{{2}}}^{2}}}+\frac{1}{6}\,{\frac {{\pi }^{2}}{{l_{{2}}}^{2}}}
\]
\\
Das gleiche Prozedere jetzt mit der Dreifach-Summe
\[
\sum _{\ell_{{3}}=3}^{\infty} \left( \sum _{\ell_{{2}}=2}^{\ell_{{3}}-1} \left( \sum _
{\ell_{{1}}=1}^{\ell_{{2}-1}} \frac{1}{\ell_1^2 \cdot\ell_2^2 \cdot\ell_3^2}  \right)  \right) 
\]
und dem dazugehörigen MAPLE Code
\begin{verbatim}
sum(sum(sum(1/(l[1]^2*l[2]^2*l[3]^2),
l[1]=1..l[2]-1),l[2]=2..l[3]-1),l[3]=3..infinity);
\end{verbatim}
Es führt zu einem ähnlichen Ergebnis. Wir erhalten nicht  $\frac{\pi^6}{7!}$, sondern den Term
\[
\sum _{l_{{3}}=3}^{\infty } \left( -\frac{1}{6}\,{\frac {{\pi }^{2}\Psi
 \left( 1,l_{{3}} \right) }{{l_{{3}}}^{2}}}+\frac{1}{6}\,{\frac {{\pi }^{2}
 \left( -1+1/6\,{\pi }^{2} \right) }{{l_{{3}}}^{2}}}+\sum _{l_{{2}}=2
}^{l_{{3}}-1}-{\frac {\Psi \left( 1,l_{{2}} \right) }{{l_{{2}}}^{2}{l_
{{3}}}^{2}}} \right)
\]
\\
Außer im Falle $\frac{\pi^2}{3!}$ liefert MAPLE kein rationales Vielfaches einer  $\pi$-Potenz. Wir erhalten stattdessen immer kompliziertere Terme mit der $\Psi$-Funktion.
$ \\ $
$ \\ $
Nun überprüfen wir, wie MAPLE mit der numerischen Berechnung der unendlichen Summen zurechtkommt. In MAPLE wird die Anzahl der Dezimalstellen mit denen gerechnet wird durch
\begin{verbatim} Digits
\end{verbatim}
festgelegt. Wir wollen mit 50 Dezimalstellen rechnen und setzen deshalb
\begin{verbatim}
Digits:=50;
\end{verbatim}
MAPLE führt dann die Berechnungen mit 50 Stellen Genauigkeit aus. Wir benötigen noch
\begin{verbatim}
evalf() 
\end{verbatim}
damit der berechnete Wert mit dieser Stellenzahl angezeigt wird. 
$ \\ $
$ \\ $
Als erstes lassen wir uns $\frac{\pi^4}{5!}$ und den Wert der dazugehörigen Summe berechnen
\begin{verbatim}
Digits:=50;
evalf(Pi^4/5!);
evalf(sum(sum(1/(l[1]^2*l[2]^2),
l[1]=1..l[2]-1),l[2]=2..infinity));
\end{verbatim}
$ \\ $
MAPLE gibt aus:
$ \\ $
\centerline{$Digits:=50$ } 
\centerline{$0.81174242528335364363700277240587592708106321393904$}
\centerline{$0.81174242528335364363700277240587592708106321393905$}
$ \\ $
Die beiden Werte sind sehr exakt. Sie unterscheiden sich nur in der letzten Ziffer.
$ \\ $
$ \\ $
Als nächstes berechnen wir  $\frac{\pi^6}{7!}$ und die dreifach ineinander geschachtelte Summe
\begin{verbatim}
evalf(Pi^6/7!);
evalf(sum(sum(sum(1/(l[1]^2*l[2]^2*l[3]^2),
l[1]=1..l[2]-1),l[2]=2..l[3]-1),l[3]=3..infinity));
\end{verbatim}
Wir erhalten:
\\
\centerline{$0.19075182412208421369647211183579759898159077938116$}
\centerline{$0.19075182412208421369647211183579759898159077938116$}
$ \\ $
Die beiden Werte sind identisch. Es stimmen alle Ziffern überein
$ \\ $
$ \\ $
Wir führen noch eine Berechnung mit  $\frac{\pi^8}{9!}$  und der dazugehörigen Vierfachsumme aus
\begin{verbatim}
evalf(Pi^8/9!);
evalf(sum(sum(sum(sum(1/(l[1]^2*l[2]^2*l[3]^2*l[4]^2),
l[1]=1..l[2]-1),l[2]=2..l[3]-1),l[3]=3..l[4]-1),l[4]=4..infinity));
\end{verbatim}
MAPLE liefert das Ergebnis:
$ \\ $
$ \\ $
\centerline{$0.026147847817654800504653261419496157949452103923173$}
\centerline{$0.026147847817654800504653261419496157949452103923173$}
$ \\ $
Auch hier sind beide Werte gleich. Es ist erstaunlich, dass MAPLE unendliche Summen mit so hoher Genauigkeit berechnen kann.
\section*{6\quad Appendix}
Nun kommen wir zum angekündigten Beweis aus Abschnitt 2. Wir wollen beweisen, dass sich das Produkt $\prod_{k=1}^{M}(1+x_{k}t)$ durch das folgende auf der rechten Seite angegeben Polynom darstellen lässt, dass also gilt:
\begin{align*}
\prod_{k=1}^{M}(1+x_{k}t)=1
+&\left[\sum_{\ell_1=1}^{M}x{_{\ell_1}}\right] t
\\
+&\left[\sum _{\ell_{{2}}=2}^{M} \left( \sum _{\ell_{{1}}=1}^{\ell_{{2}}-1}x_{{\ell_{{1}}}
} x_{{\ell_{{2}}}} \right)\right] t^2
\\
+&\left[\sum _{\ell_{{3}}=3}^{M} \left( \sum _{\ell_{{2}}=2}^{\ell_{{3}}-1} \left( 
\sum _{\ell_{{1}}=1}^{\ell_{{2}}-1}x_{{\ell_{{1}}}} x_{{\ell_{{2}}}} x_{{\ell_{{3}}}} \right)  \right)\right] t^3
\\
+&\left[{\sum _{\ell_{{4}}=4}^{M} \left(\sum _{\ell_{{3}}=3}^{\ell_4-1} \left( \sum _{\ell_{{2}}=2}^{\ell_{{3}}-1} \left( 
\sum _{\ell_{{1}}=1}^{\ell_{{2}}-1} x_{{\ell_{{1}}}}  x_{{\ell_{{2}}}} x_{{\ell_3}} x_{{\ell_4}} \right)  \right)\right)}\right] t^4
\\
&\qquad\qquad\qquad\qquad{\vdots}
\\
+&\left[\sum\limits_{{\ell_M} = M}^M {\left( {\sum\limits_{{\ell_{M - 
1}} = M - 1}^{{\ell_M} - 1}  \cdots \left( {\sum\limits_{{\ell_2} = 2}
^{{\ell_3} - 1} {\left( {\sum\limits_{{\ell_1} = 1}^{{\ell_{2} - 1}} 
x_{\ell_1} x_{\ell_2} \cdots x_{\ell_{M-1}} x_{\ell_M} } \right)} } 
\right) \cdots } \right)}\right]t^M
\end{align*}

Der Beweis erfolgt durch Induktion nach $M$. Um die aufwendige 
Schreibarbeit für die ineinander geschachtelten Summen zu vermeiden,
 verwenden wir die für elementar-symmetrische Polynome übliche 
$\sigma$-Notation:
\[
\sum_{\ell_1=1}^{M}x{_{\ell_1}}
=\sigma_{M,1}
\]
\[
\sum _{\ell_{2}=2}^M 
\left(
\sum _{\ell_1=1}^{\ell_2-1}
x_{\ell_1} x_{\ell_2}
\right)
=\sigma_{M,2}
\]
\[
\sum _{\ell_{3}=3}^M
\left(\sum _{\ell_2=2}^{\ell_3-1}
\left(\sum _{\ell_1=1}^{\ell_2-1}
x_{\ell_1}x_{\ell_2}x_{\ell_3}
\right)\right)
=\sigma_{M,3}
\]
\[\vdots\]
\[
\sum\limits_{{\ell_k} = k}^M {\left( {\sum\limits_{{\ell_{k - 1}} = k - 1}^{{\ell_k} - 1}  \cdots  \left( {\sum\limits_{{\ell_2} = 2}^{{\ell_3} - 1} {\left( {\sum\limits_{{\ell_1} = 1}^{{\ell_{2} - 1}} x_{\ell_1}\cdot x_{\ell_2} \cdots x_{\ell_{k-1}} \cdot x_{\ell_k} } \right)} } \right) \cdots } \right)}
=\sigma_{M,k}
\]
\[\vdots\]
\[
\sum _{\ell_{M}=M}^M
\left(\sum _{\ell_{M-1}=M-1}^{\ell_M-1}
\ldots
\left(\sum _{\ell_2=2}^{\ell_3-1}
\left( \sum _{\ell_1=1}^{\ell_2-1}
x_{\ell_1}x_{\ell_2}\ldots x_{\ell_M-1}x_{\ell_M} \right)\right)
\ldots\right)
=\sigma_{M,M}
\]
\[\]
Dadurch erhalten wir die einfache Summe
\begin{align*}
\prod_{k=1}^{M}(1+x_{k}t)=1
+\sigma_{M,1}\; t
+\sigma_{M,2}\; t^2
+\ldots
+\sigma_{M,M-1}\; t^{M-1}
+\sigma_{M,M}\; t^M=1+\sum_{k=1}^{M}\sigma_{M,k}\; t^{k}
\end{align*}
$ \\ $
$ \\ $
\textbf{Beweis}
$ \\ $
\underline{Step 1:}
\\
Es sei $M = 1$. Auf der linken Seite steht dann 
\[\prod_{k=1}^{1}(1+x_{1}t)=1+x_{1}t\]
und auf der rechten Seite erhalten wir
\[
1+\sigma_{1,1}t=1+\left[\sum_{\ell_1=1}^{1}x{_{\ell_1}}\right]t=1+x_{1}t.
\]
Beide Seiten sind gleich, somit ist der Satz richtig für $M = 1$.
$ \\ $
$ \\ $
\underline{Step 2:}
$ \\ $
Wir nehmen an, die Aussage ist richtig für $M = m \geq 1$, d.h. es gilt
\begin{align*}
\prod_{k=1}^{m}(1+x_{k}t)
=1+\sum_{k=1}^{m}\sigma_{m,k}\; t^{k}
\end{align*}
$ \\ $
\underline{Step 3:}
$ \\ $
Wir wollen beweisen, dass die Aussage auch für $m+1$ richtig ist, dass also gilt
\begin{align*}
\prod_{k=1}^{m+1}(1+x_{k}t)
=1+\sum_{k=1}^{m+1}\sigma_{m+1,k}\; t^{k}
\end{align*}
$ \\ $
Wie kommen wir von Step $2\,$ zu Step $3\,$? Sehen wir uns einmal die linke Seit an. Wir stellen fest, dass sie sich wie folgt schreiben lässt.
\[
\prod_{k=1}^{m+1}\left(1+x_{k}t\right)
=\left(\prod_{k=1}^m\left(1+x_{k}t\right)\right)\left(1+x_{m+1}t\right)
\]
Für das Produkt
\begin{align*}
\prod_{k=1}^{m}(1+x_{k}t)
\end{align*}
können wir die Induktionsannahme
\begin{align*}
1+\sum_{k=1}^{m+1}\sigma_{m+1,k}\; t^{k}
\end{align*}
einsetzen und erhalten
\begin{align*}
\left(1+\sum_{k=1}^{m}\sigma_{m,k}\; t^{k}\right)
\left( 1+x_{m+1}t\right)
\end{align*}
$ \\ $
Wir lösen die Klammer auf und fassen so weit wie möglich zusammen
$ $
\begin{align*}
\left(1+\sum_{k=1}^{m}\sigma_{m,k}\; t^{k}\right)
\left( 1+x_{m+1}t\right)
&=1+\sum_{k=1}^{m}\sigma_{m,k}\; t^{k}
+x_{m+1}t+\left(\sum_{k=1}^{m}\sigma_{m,k}\; t^{k}\right)x_{m+1}t
\\
&=1+x_{m+1}t+\sum_{k=1}^{m}\sigma_{m,k}\; t^{k}
+\sum_{k=1}^{m}\sigma_{m,k}\; x_{m+1} \; t^{k+1}
\end{align*}
Um Schreibarbeit zu ersparen kürzen wir die Terme auf der rechten Seite wie folgt ab
\[
A=x_{m+1}t+\sum_{k=1}^{m}\sigma_{m,k}\; t^{k}
\]
\[
B=\sum_{k=1}^{m}\sigma_{m,k}\; x_{m+1} \; t^{k+1};
\]
Die rechte Seite schreibt sich dann einfach als
\[
1+A+B
\]
Als erstes zerlegen wir den Term $A$ in zwei Summen.
\begin{align*}
A=x_{m+1}t+\sum_{k=1}^{m}\sigma_{m,k}\; t^{k}
&=x_{m+1}t+\sum_{k=1}^{1}\sigma_{m,k}\; t^{k}
+\sum_{k=2}^{m}\sigma_{m,k}\; t^{k}
\end{align*}
und fassen zusammen
\begin{align*}
A
=x_{m+1}t+\sigma_{m,1}\;t+\sum_{k=2}^{m}\sigma_{m,k}\;t^{k}
=\left(x_{m+1}+\sigma_{m,1}\right)t+\sum_{k=2}^{m}\sigma_{m,k}\;t^{k}\end{align*}
\\
Nun gilt
\[
\sigma_{m,1}
=\sum_{\ell_1=1}^{m}x_{\ell_1}
\]
Dies eingesetzt ergibt
\begin{align*}
A
=\left(x_{m+1}+\sum_{\ell_1=1}^{m}x_{\ell_1}\right)t+\sum_{k=2}^{m}\sigma_{m,k}\;t^{k}
=\sum_{\ell_1=1}^{m+1}x_{\ell_1}t+\sum_{k=2}^{m}\sigma_{m,k}\;t^{k}
\end{align*}
Jetzt nehmen wir uns den zweiten Term vor. Es war
\[
B=\sum_{k=1}^{m}\sigma_{m,k}\; x_{m+1} \; t^{k+1}
\]
Wir zerlegen den Term in zwei Summen und nehmen eine Indexverschiebung vor.
\begin{align*}
B
&=\sum _{k=1}^{m}\sigma_{{m,k}}\;x_{{m+1}}\;{t}^{k+1}
=\sum _{k=1}^{m-1}\sigma_{{m,k}}\;x_{{m+1}}\;{t}^{k+1}
+\sum _{k=m}^{m}\sigma_{{m,k}}\;x_{{m+1}}\;{t}^{k+1}\\
&=\sum _{k=1}^{m-1}\sigma_{{m,k}}\;x_{{m+1}}{t}^{k+1}
+\sigma_{m,m} \cdot x_{m+1}\;t^{m+1}\\
&=\sum _{k=2}^{m}\sigma_{{m,k-1}}\;x_{{m+1}}\;{t}^{k}
+\left(\prod_{k=1}^{m}x_{k}\right) x_{m+1}t^{m+1}\\
&=\sum _{k=2}^{m}\sigma_{{m,k-1}}\;x_{{m+1}}\;{t}^{k}
+\left(\prod_{k=1}^{m+1}\;x_{k}\right)t^{m+1}\\
&=\sum _{k=2}^{m}\sigma_{{m,k-1}}\;x_{{m+1}}\;{t}^{k}
+\sum _{k=m+1}^{m+1}\sigma_{{m,k}}\;{t}^{k+1}
\end{align*}
\\
Bis jetzt haben wir folgendes erreicht
\begin{align*}
&A=\sum_{\ell_1=1}^{m+1}x_{\ell_1}t+\sum_{k=2}^{m}\sigma_{m,k}\;t^{k}\;t^{k}\\
&B=\sum _{k=2}^{m}\sigma_{{m,k-1}}\;x_{{m+1}}\;{t}^{k}
+\sum _{k=m+1}^{m+1}\sigma_{{m,k}}\;{t}^{k+1}
\end{align*}
Eingesetzt erhalten wir
\begin{align*}
1+A+B
&=1+\sum _{k=1}^{1}\sigma_{{m+1,k}}\;{t}^{k}
+\sum _{k=2}^{m}\sigma_{{m,k}}\;{t}^{k}
+\sum _{k=2}^{m}\sigma_{{m,k-1}}x_{{m+1}}\;{t}^{k}
+\sum _{k=m+1}^{m+1}\sigma_{{m,k}}\;{t}^{k+1}\\
&=1+\sum _{k=1}^{1}\sigma_{{m+1,k}}\;{t}^{k}
+\sum _{k=2}^{m}\left(\sigma_{{m,k}}+\sigma_{{m,k-1}}x_{{m+1}}\right)\;{t}^{k}
+\sum _{k=m+1}^{m+1}\sigma_{{m,k}}\;{t}^{k+1}\end{align*}
\\
Auf die zweite Summe wenden wir Lemma 1 an:
\[
\sum _{k=2}^{m}\left(\sigma_{{m,k}}+\sigma_{{m,k-1}}x_{{m+1}}\right)\;{t}^{k}=\sum_{k=2}^{m}\sigma_{m+1,k}\; t^{k}
\]
Setzen wir dies ein, so ergibt sich
\[
1+\sum_{k=1}^{1}\sigma_{m+1,k}\; t^{k}
+\sum_{k=2}^{m}\sigma_{m+1,k}\; t^{k}
+\sum_{k=m+1}^{m+1}\sigma_{m+1,k}\; t^{k}
\]
Die drei Summen können wir zu einer einzigen Summe zusammenfassen:
\[
\left(1+\sum_{k=1}^{m}\sigma_{m,k}\; t^{k}\right)
\left( 1+x_{m+1}t\right)=1+\sum_{k=1}^{m+1}\sigma_{m+1,k}\; t^{k}
=\prod_{k=1}^{m+1}(1+x_{k}t)\]
Das ist genau das, was wir beweisen wollten. Somit ist die Aussage für alle $M\in\mathbb N$ richtig.

\end{document}